\documentclass{amsart}
\usepackage{t1enc}
\usepackage[utf8]{inputenc}
\usepackage{graphicx}

\usepackage{amsmath}
\usepackage{amssymb}
\usepackage{amsthm}
\usepackage{mathrsfs}
\usepackage[all]{xy}
\usepackage{enumerate}

\renewcommand{\epsilon}{\varepsilon}
\theoremstyle{plain}
\newtheorem{thm}{Theorem}[section]
\newtheorem{prop}[thm]{Proposition}

\newtheorem{cor}[thm]{Corollary}
\newtheorem{lem}[thm]{Lemma}

\theoremstyle{definition}
\newtheorem{dfn}[thm]{Definition}

\theoremstyle{remark}
\newtheorem*{rem}{Remark}

\def\doublear[#1]{\ar@<2pt>[#1]\ar@<-2pt>[#1]}

\DeclareMathOperator{\Top}{Top}
\DeclareMathOperator{\pt}{pt}

\DeclareMathOperator{\Aut}{Aut}

\DeclareMathOperator{\Gal}{Gal}

\DeclareMathOperator{\Spec}{Spec}

\DeclareMathOperator{\ga}{\pi_1^{alg}}
\DeclareMathOperator{\gtop}{\pi_1^{top}}

\DeclareMathOperator{\gtemp}{\pi_1^{temp}}
\DeclareMathOperator{\gtempL}{\pi_1^{\mbb L-temp}}
\DeclareMathOperator{\glog}{\pi_1^{log}}
\DeclareMathOperator{\ggeom}{\pi_1^{log-geom}}

\DeclareMathOperator{\KCov}{KCov}
\DeclareMathOperator{\KCovgeom}{KCov_{geom}}

\DeclareMathOperator{\Covalg}{Cov^{alg}}
\DeclareMathOperator{\Covtop}{Cov^{top}}
\DeclareMathOperator{\Covtemp}{Cov^{temp}}
\DeclareMathOperator{\CovtempL}{Cov^{\mbb L-temp}}
\DeclareMathOperator{\ket}{k\acute{e}t}

\DeclareMathOperator{\rk}{rk}
\DeclareMathOperator{\gp}{gp}

\DeclareMathOperator{\Coker}{Coker}

\DeclareMathOperator{\DDtemp}{DD_{temp}}
\DeclareMathOperator{\Dtop}{\mcal D_{top}}

\DeclareMathOperator{\Dtemp}{\mcal D_{temp}}
\DeclareMathOperator{\DtempL}{\mcal D_{\mbb L-temp}}

\DeclareMathOperator{\Hom}{Hom}

\DeclareMathOperator{\Ens}{Set}
\DeclareMathOperator{\Set}{Set}
\DeclareMathOperator{\tSet}{\text{-}Set}

\DeclareMathOperator{\tfEns}{\text{-}fSet}

\DeclareMathOperator{\fSet}{fSet}

\DeclareMathOperator{\Graph}{Graph}
\DeclareMathOperator{\GenGraph}{GenGraph}
\DeclareMathOperator{\Stab}{Stab}

\DeclareMathOperator{\an}{an}

\DeclareMathOperator{\op}{op}

\newcommand\projLim{\mathop{\mathrm{\underset{\longleftarrow}{Lim}}}}
\newcommand\injLim{\mathop{\mathrm{\underset{\longrightarrow}{Lim}}}}

\newcommand{\mring}{\mathring}
\DeclareMathOperator{\triv}{tr}
\newcommand{\tr}{\triv}

\DeclareMathOperator{\OutGptop}{OutGp_{top}}
\DeclareMathOperator{\Pt}{Pt}

\newcommand{\ie}{\emph{i.e.} }
\newcommand{\findem}{\end{proof}}
\newcommand{\dem}{\begin{proof}}

\newcommand{\da}{\begin{displaystyle}}
\newcommand{\db}{\end{displaystyle}}
\newcommand{\dar}{\downarrow}

\newcommand{\mcal}{\mathcal}
\newcommand{\mbf}{\mathbf}

\newcommand{\mbb}{\mathbb}
\newcommand{\fk}{\mathfrak}

\DeclareMathOperator{\C}{C}

\newcommand*{\cube}[8]{\xymatrix {
    #1 \ar[rr] \ar[dd] \ar[dr] && #2 \ar[dr] \ar[dd]|\hole \\
    & #3 \ar[rr] \ar[dd] && #4 \ar[dd] \\
    #5 \ar[rr]|\hole \ar[dr] && #6 \ar[rd] \\
    & #7 \ar[rr] && #8 \\
  }}

\begin{document}

\title[Cospecialization of tempered groups for curves]
{Coverings in $p$-adic analytic geometry and log covers I:\\
Cospecialization of the $(p')$-tempered fundamental group for a family of curves}
\author{Emmanuel Lepage}
\email{lepage@math.jussieu.fr}
\address{Université Pierre et Marie Curie\\
Institut Mathématique de Jussieu\\
4 place Jussieu\\
75005 PARIS, FRANCE}
\begin{abstract}
The tempered fundamental group of a $p$-adic analytic space classifies
covers that are dominated by a topological cover (for the Berkovich
topology) of a finite \'etale cover of the space. Here we construct
cospecialization homomorphisms between $(p')$ versions of the tempered fundamental
groups of the fibers of a smooth family of curves with semistable reduction. To
do so, we will translate our problem in terms of cospecialization morphisms of
fundamental groups of the log fibers of the log reduction and we
will prove the invariance of the geometric log fundamental group of
log smooth log schemes over a log point by change of log
point.
\end{abstract}
\maketitle
\section*{Introduction}
 In general topology, the fundamental group  
of a connected locally contractible pointed space classifies its (unramified)  
covers. A. Grothendieck developed an avatar in abstract algebraic  
geometry: he attached to any algebraic variety a profinite  
fundamental group, which classifies its finite \'etale covers. For
a complex algebraic variety, Grothendieck's fundamental group is
canonically isomorphic to the profinite completion of the topologic
fundamental group of the corresponding topological space.

Here we are interested in an analog in $p$-adic geometry. More precisely we
will study the \emph{tempered fundamental group} of $p$-adic varieties
defined by Y. André.  The profinite completion of the tempered fundamental
group of any smooth $p$-adic algebraic variety coincides with Grothendieck's algebraic
fundamental group. It also accounts for the usual (infinite) ``uniformizations'' in
$p$-adic analytic geometry such as the uniformization of Tate elliptic curves,
which are historically at the very basis of $p$-adic rigid
geometry. Such uniformizations give infinite discrete quotients
of the tempered fundamental group.

The framework of this paper for non-archimedean analytic geometry will be Berkovich spaces. The underlying space of the analytification of an affine algebraic variety $\Spec A$ in the sense of Berkovich is the set of multiplicative seminorms on $A$ with value in $\mbf R_{\geq 0}$ extending the norm of the base field, endowed with the coarsest topology that makes the evaluation of the norm of any element $f\in A$ continuous. The analytification of a smooth algebraic variety is locally contractible, which ensures the existence of universal topological covers. 
Since the analytification (in the sense of V. Berkovich or of rigid geometry)
of a finite étale cover of a $p$-adic algebraic variety is not
necessarily a topological cover, André had to consider a category of
covers slightly bigger than just the category of topological covers.
He defined tempered covers, which are (possibly infinite) \'etale covers in
the sense of A.J. de Jong (that is to say, which are,  
Berkovich-locally on the base, direct sums of finite covers)
 such that, after pulling
back by some finite \'etale cover, they become topological covers (for the
Berkovich topology). The \emph{tempered fundamental group} is the prodiscrete group that classifies those tempered covers.  To give a more handful description, if
one has a sequence of pointed finite Galois connected covers $((S_i,s_i))_{i\in \mbf N}$ such that the corresponding pointed pro-cover of $(X,x)$ is the
universal pro-cover of $(X,x)$, and if $(S^{\infty}_i,s^{\infty}_i)$ is a
universal topological cover of $S_i$, the tempered fundamental group of
$X$ can be seen as $\gtemp(X,x)=\varprojlim_i \Gal(S^{\infty}_i/X)$.
Therefore, to understand the tempered fundamental group of a variety, one
has to understand the topological behavior of its finite \'etale
covers.

In the case of a curve, the question becomes more concrete since there
is a natural embedding of the geometric realization of the graph of its
stable model into the Berkovich space of the curve which is a homotopy
equivalence.

Among applications of tempered fundamental groups, let us cite in passing
the theory of $p$-adic orbifolds and $p$-adic triangle groups \cite{andre1}
and a
$p$-adic version of Grothendieck-Teichm\"uller theory \cite{andre2}.

\bigskip

In this article we will be interested in the variation of the tempered
fundamental group of the fibers of a family of curves. This article will be
followed by another one~\cite{cosp2}, in which we will consider higher dimensional
families.

If $\bar y_1,\bar y_2$ are to geometric points of a scheme $Y$, a specialization $\bar y_2\to\bar y_1$ is a $Y$ morphism from $\bar y_2$ to the strict localization $Y(\bar y_1)$ of $Y$ at $\bar y_1$. Equivalently, a specialization $\bar y_2\to\bar y_1$ is a morphism of functors $(\ )_{\bar y_1}\to (\ )_{\bar y_2}$, where $(\ )_{\bar y}$ is the functor from the \'etale topos of $Y$ to the category of sets that maps an \'etale sheaf $\mcal F$ to its stalk $\mcal F_{\bar y}$.
For a proper morphism of schemes $f: X\to Y$ with geometrically connected
fibers and a specialization $ 
\bar y_2 \to \bar y_1$ of geometric points of $Y$, A. Grothendieck  
has defined algebraic fundamental groups $\ga(X_{\bar y_i})$  
and a specialization homomorphism  $\ga(X_{\bar y_1})\to  
\ga(X_{\bar y_2})$. Grothendieck's specialization theorem 
tells that this homomorphism is surjective if $f$ is separable and induces an  
isomorphism between the prime-to-$p$ quotients if $f$ is smooth  
(here, $p$ denotes the characteristic of $\bar y_2$), cf. \cite[cor. X.2.4,
cor. X.3.9]{sga}.

In complex analytic geometry, a smooth and proper morphism is locally a trivial
fibration of real differential manifolds, so that, in particular, all the
fibers are homeomorphic, and thus have isomorphic (topological) fundamental
groups.

The aim of this paper is to find some analog of the specialization
theorem of Grothendieck in the case of the tempered fundamental group.

In this paper, we will concentrate on the case of curves.

\bigskip

One problem which appears at once in looking for some non archimedean analog of
Grothendieck's specialization theorems is that there are in general no non trivial
specializations between distinct points of a non archimedean analytic
(Berkovich or rigid) space: for example a separated Berkovich space has a
Hausdorff underlying topological space, so that if there is a
cospecialization (for the Berkovich topology, the étale topology\dots)
between two geometric points of a Berkovich space, the two geometric points
must have the same underlying point. Thus we will assume we have a
model over the ring of integers of our non-archimedean field (with good enough properties) and we will look at
the specializations in the special fiber.

\smallskip
 
We want to understand how the tempered fundamental group of the geometric fibers of a
smooth family varies. Let us for instance consider a family of elliptic curves. The tempered
fundamental group of an elliptic curve over a complete algebraically closed
non archimedean closed field is $\widehat{\mbf Z}^2$ if it has
good reduction, and $\widehat{\mbf Z}\times \mbf Z$ if it is a Tate
curve. In particular, by looking at a moduli space of stable pointed elliptic curves
with level structure\footnote{to avoid stacks. However, the cospecialization
  homomorphisms we will construct will be local for the étale topology of
  the special fiber of the base. Thus, the fact that the base is a
  Deligne-Mumford stack is not really a problem.}, the tempered fundamental group (or any reasonable
$(p')$-version) cannot be constant.

Moreover, if one looks at the moduli space over $\mbf Z_p$, and 
considers a curve $E_1$ with bad reduction and a curve $E_0$ with generic
reduction (hence good reduction), there cannot be a morphism
$\gtemp(E_0)\to\gtemp(E_1)$ which induces Grothendieck's specialization on
the profinite completion, although the reduction point corresponding to
$E_1$ specializes to the reduction point corresponding to $E_0$. Therefore
there cannot be any reasonable specialization theory.

On the contrary, if
one has two geometric points $\eta_1$ and $\eta_2$ of the moduli space such
that the reduction of $\eta_1$ specializes to the reduction of $\eta_2$,
then $E_{\eta_1}$ has necessarily better reduction than $E_{\eta_2}$ and
there is some morphism
$\gtemp(E_{\eta_2})\to\gtemp(E_{\eta_1})$ that induces an isomorphism between
the profinite completions. Thus we want to look for a \emph{cospecialization} of the tempered fundamental group.

\smallskip
 
The topological behavior of general finite étale covers is too
complicated to hope to have a simple cospecialization theory without
adding a $(p')$ condition on the covers: for example
two Mumford curves over some finite extension of $\mbf Q_p$ with isomorphic
geometric tempered fundamental group have the same metrized graph of stable reduction \cite{metricmumford}. Thus even if two Mumford curves
have isomorphic stable reduction (and thus the point corresponding to their stable reduction is the same), they may not have isomorphic
tempered fundamental group in general. Thus we will only
study here finite covers that are dominated by a finite Galois cover
whose order is prime to $p$, where $p$ is the residual characteristic
(which can be $0$; such
a cover will be called a $(p')$-finite cover). Then, it becomes
natural to introduce a $(p')$-tempered fundamental group which classifies
tempered covers that become topological covers after pullback along
some $(p')$-finite cover. It should be remarked that this
$(p')$-tempered fundamental group cannot be in general recovered from the
tempered fundamental group.

The $(p')$-tempered fundamental group of a curve was already studied by
S. Mochizuki in~\cite{mochi}. It can be described in terms of a graph of profinite
groups. From this description, one easily sees that the isomorphism class
of the $(p')$-tempered
fundamental group of a $p$-adic curve depends only of the stratum of the
Knudsen stratification of the moduli space of stable curves in which the
stable reduction lies. Moreover if one has two strata $x_1$ and $x_2$ in the moduli space
of stable curves such that $x_1$ is in the closure of $x_2$, one can easily
construct  morphisms from the graph of groups
corresponding to $x_1$ to the graph of groups corresponding to $x_2$
(inducing morphisms of tempered fundamental groups which induce 
isomorphisms of the pro-$(p')$ completions).

\subsection*{}
We shall study the following situation. Let $O_K$ be a complete discrete
valuation ring, $K$ be its fraction field, $k$ be its residue field and $p$ be
its characteristic (which can be $0$). A proper semistable pointed curve $(X,D)$ over a scheme $S$ is given by a flat and proper morphism $X\to S$ with semistable geometric fibers, and $D$ is a closed subscheme of $X$ which lies inside the smooth locus of $X\to S$ and such that $D\to S$ is \'etale.
Let $(X,D)$ be a proper semistable pointed curve over $O_K$ smooth over $K$ and let $U=
X_{\eta}\backslash D_{\eta}$. Let us describe the tempered
fundamental group of $U^{\an}_{\bar \eta}$ in terms of $X_{\bar s}$ \cite{mochi}.

 Let us make sure at first that we can get such a
description for the pro-$(p')$ completion, \ie the algebraic
fundamental group. One cannot apply directly Grothendieck's specialization
theorems (even if $U=X_{\eta}$) since the special fiber is not smooth but only
semistable. Indeed, a pro-$(p')$ geometric cover of the generic fiber will
generally only induce a Kummer cover on the special fiber. These are
naturally described in terms of \emph{log geometry}, more precisely in
terms of \emph{Kummer-étale} (két) covers of a log scheme. One can endow $X$ (and thus
$X_s$ too by restriction)
with a natural log structure such that the pro-$(p')$ fundamental group of
$U$ is isomorphic to a pro-$(p')$ log fundamental group (as defined
in~\cite{ill}) of $X_{\bar s}$. One then gets a description of
$\ga(U_{\bar\eta})$ by taking the projective limit under tame covers of
$K$, or equivalently under két covers of $s$ endowed with its natural
log structure: there is an equivalence between finite étale covers of $U_{\bar
  \eta}$ and ``geometric két covers'' of $X_{\bar s}$.

A két cover of $X$ will still be a semistable model of its generic fiber
if one replaces $K$ by some tame extension. Thus, one can describe the
topology of the corresponding cover of $U_{\bar \eta}$ in terms of the
graph of the corresponding geometric két cover of $X_{\bar s}$.

\smallskip

Let us now come back to the problem of cospecialization. Let $X\to Y$ be a semistable
curve over $O_K$ with $X\to Y$ endowed with compatible log structure (see definition~\ref{stlogcurve}).

Let $\bar \eta_1$ (resp. $\bar \eta_2$) be a (Berkovich) geometric point of
$Y_0:=Y^{\an}_{\tr}\cap \fk Y_\eta\subset Y^{\an}_K$, where $Y_{\tr}$ is the locus of $Y$ where the log structure is trivial and $\fk Y_\eta$ is the generic fiber of the formal completion of $Y$ along its closed fiber (if $Y$ is proper, then $Y_0=Y^{\an}_{\tr}$). Let $\bar s_1$ (resp. $\bar s_2$) be its log reduction in $Y_s$ .

To use the previous description of the tempered fundamental group of
$U_{\bar \eta_1}$ and $U_{\bar \eta_2}$ in terms of $X_{\bar s_1}$ and
$X_{\bar s_2}$, we have to assume that $\bar \eta_1$ and $\bar \eta_2$ lie
over Berkovich points with discrete valuation.

\bigskip

The main result of this paper is the following:
\begin{thm}[{th. \ref{cospcourbes}}]\label{sp}Let $K$ be a complete discretely valued field. Let $\mbb L$ be a set of primes that does not contain the residual characteristic of $K$. Let $Y\to \Spec O_K$ be a morphism of log schemes of finite type. Let $Y_0=Y^{\an}_{\tr}\cap \fk Y_\eta\subset Y^{\an}$ where $\fk Y$ is the completion of $Y$ along its closed fiber. Let $X\to Y$ be a proper semistable
curve  with compatible log structure. Let $U=X_{\tr}$. Let $\eta_1$ and $\eta_2$ be two Berkovich points of $Y_0$ whose residue fields have discrete valuation, and let $\bar \eta_1,\bar \eta_2$ be
  geometric points above them. 
Let $\bar s_2\to\bar s_1$ be a log specialization
  of their log reductions such that there exists a compatible specialization $\bar \eta_2\to\bar \eta_1$, then there is a
  cospecialization homomorphism $\gtempL(U_{\bar\eta_1})\to\gtempL(U_{\bar\eta_2})$. Moreover, it is an isomorphism if $\overline M_{Y,\bar s_1}\to\overline M_{Y,\bar s_2}$ is an isomorphism.\end{thm}

Let us come back to our example of the moduli space of pointed stable elliptic curves with
high enough level
structure $M$ over $O_K$, and let $C$ be the canonical stable elliptic curve on
$M$. Let $\mbb L$ be a set of primes that does not contain the residual characteristic of $K$. If $\eta_1$ and $\eta_2$ are two Berkovich points of $M_\eta$, they
are in $M_\eta^{\tr}$ if and only if $C_{\eta_1}$ and $\C_{\eta_2}$ are
smooth. $C\to M$, endowed with their natural log-structures over
$(O_K,O_K^*)$, is a semistable morphism of log schemes. One thus get a cospecialization
outer morphism $\gtempL(C_{\overline \eta_1})\to\gtempL(C_{\overline \eta_2})$ for
every specialization $\overline s_2\to\overline s_1$, which is an
isomorphism if $\overline s_1$ and $\overline s_2$ are in the same stratum
of $M_s$. Since the moduli stack
of pointed stable elliptic curves over $\Spec k$ has only two strata, one
corresponding to smooth elliptic curves $M_0$ and one to singular curves $M_1$,
one gets that $\gtempL(E_1)\simeq\gtempL(E_2)$ if $E_1$ and $E_2$
are two curves with good reduction or two Tate curves (the isomorphism
depends on choices of cospecializations). Since $M_1$ is in the closure of
$M_0$ one gets a morphism from the tempered fundamental group of a Tate
curve to the tempered fundamental group of an elliptic curve with good
reduction.

\smallskip

The first thing we need in order to construct the cospecialization homomorphism for
tempered fundamental groups is a
specialization morphism between the $(p')$-log geometric fundamental groups of
$X_{\bar s_1}$ and $X_{\bar s_2}$. Such a specialization morphism will be
constructed by proving that one can extend any $(p')$-log geometric
cover of $X_{s_1}$ to a két cover of $X_U$ where $U$ is some két
neighborhood of $s_1$. If one has
such a specialization morphism, by comparing it to the fundamental groups
of $X_{\bar\eta_1}$ and $X_{\bar\eta_2}$ and using Grothendieck's
specialization theorem, we will easily get that it must be an
isomorphism. This specialization morphism is easily deduced
from~\cite{org2} if $s_1$ is a strict point of $Y$ (\ie the log structure
of $s_1$ is simply the one induced by $Y$), \ie the log
structure of $s_1$ is just the pull back of the log structure of $Y$, but
is not straightforward when the log structure is really modified. Thus
we will study the invariance of the log geometric fundamental group by
change of fs base point. The main result we will prove (in any dimension) is the following~:
\begin{thm}[{th. \ref{chgmtbase1}}]\label{chbase} Let $s'\to s$ be a morphism of fs log points
  with isomorphic algebraically closed underlying fields. Let $X\to s$ be a
  saturated morphism of log schemes with $X$ noetherian and let $X'\to
  s'$ be the pull back to $s'$. Then the map
  $\ggeom(X'/s',\bar x')\to\ggeom(X/s,\bar x)$  is an isomorphism.\end{thm}
It is interesting to notice that, in this situation, this is an isomorphism
for the full fundamental group, and not only of the pro-$(p')$ part. This
mainly comes from the fact that the morphism of underlying schemes $\mring
X'\to\mring X$ is an isomorphism (so that the problem only comes from the
logarithmic structure and not the schematic structure). This result is
proved by a local study on $X$ for the strict étale topology.

\smallskip

 Then we have to construct cospecialization topological morphisms for a
 semistable curve, more precisely cospecialization morphisms of the graphs
 of the geometric fibers. This will be done étale locally. These morphisms are not morphisms of graphs
 in the usual sense, since an edge can be contracted over a vertex, but
 still give a map between their geometric realizations, whence a 
 map of homotopy types $U_{\bar\eta_1}\to U_{\bar\eta_2}$. This can also be
 done for any két cover of $X_{\bar s_1}$: we thus get such a map of homotopy
 types for every $(p')$-cover of $U_{\bar\eta_1}$. Those maps are
 compatible, and thus glue together to give the wanted cospecialization of tempered
 fundamental groups.

\subsection*{}
The paper is organized as follows.

\smallskip

In the first section, we will study specialization of fundamental groups.

In the second section, we will construct cospecialization maps of graphs of the geometric fibers of a semistable curve.

In the last section, we will prove theorem \ref{sp}.

\bigskip

This work is part of a PhD thesis. I would like to thank my advisor, Yves André, for
suggesting me to work on the cospecialization of the tempered fundamental
group and taking the time of reading and correcting this work. I would also
like to thank Luc Illusie and Fumiharu Kato for taking interest in my
problem about the invariance of geometric log fundamental groups by base change.

\section{Specialisation of log fundamental groups}\label{loggeometry}
The main result of this part will be the invariance of the log geometric
fundamental group announced in theorem~\ref{chbase}. We will deduce from it
morphisms of specialization for the pro-$(p')$ log geometric fundamental
group of the fibers of a proper log smooth saturated morphism. 

\subsection{Log fundamental groups}
For a curve with bad reduction, one cannot apply Grothendieck's
specialization theorem to describe the geometric fundamental group of the
curve in terms of the fundamental group of its stable reduction since the
family is not smooth. However, such a comparison result exists in the realm
of log geometry. More precisely, if one considers a smooth and proper variety with
semistable reduction, the semistable model can naturally be endowed with a
log structure, and the pro-$(p')$ fundamental group of the variety is
canonically isomorphic to the pro-$(p')$ log fundamental group of the
semistable reduction. Here we recall the basic definitions and results
about log fundamental groups.

\smallskip

First, recall some usual notations about monoids and log schemes. All the monoids we consider are commutative with unit.
If $P$ is a monoid, then $P^*$ is the group of invertible elements of $P$ and $P^{\gp}$ is the universal group together with a morphism of monoids $P\to P^{\gp}$. A monoid is integral if $P\to P^{\gp}$ is injective. A monoid $P$ is \emph{sharp} if $P^*$ is trivial. The sharpification $P/P^*$ is denoted by $\overline P$.
If $X$ is a log
scheme, the sheaf of monoids defining its log structure will usually be
denoted by $M_X$, the sharpification $M_X/O^*_X$ of $M_X$ will be denoted by $\overline M_X$,  the underlying scheme will be denoted by $\mring{X}$ and
the open subset of $\mring X$ where the log structure is trivial will be
denoted $X_{\tr}$.

A morphism $X\to Y$ of log schemes is \emph{strict} if the log structure on $X$ is the pullback  log structure of the log structure on $Y$. If $M_X$ and $M_Y$ are integral, then $f:X\to Y$ is strict if and only if, for every geometric point $\bar x$ of $X$, $\overline M_{Y,f(\bar x)}\to\overline M_{X,\bar x}$ is an isomorphism.

If $P$ is a monoid, one denotes by $\Spec P$ the set of primes of $P$. There is a natural map $\Spec \mbf Z[P]\to\Spec P$.

A monoid $P$ \emph{fine and saturated} (or \emph{fs} for short) if it is finitely generated, integral and, for every $a\in P^{\gp}$ such that there exists a positive integer $n$ such that $a^n\in P$, then $a\in P$.
A log scheme $X$ is \emph{fs} if, locally for the \'etale topology of $X$, there is an fs monoid $P$ and a morphism $P\to M_X$ such that $P^a\to M_X$, where $P^a$ is the log structure associated to $P\to O_X$, is an isomorphism ($P\to M_X$ is then called a \emph{fs chart} modeled on $P$). Giving a chart $P\to M_X$ is equivalent to giving a strict morphism of log schemes $X\to \Spec \mbf Z[P]$.

If $K$ is a complete discretely valued field, $S=\Spec O_K$ will be endowed in this paper with the log structure associated to $O_K\backslash\{0\}\to O_K$. If $\pi$ is a uniformizer of $O_K$, the map $\alpha:\mbf N\to O_K$ defined by $\alpha(n)=\pi^n$ is an fs chart.

If $\mbb L$ is a set of prime numbers, a $\mbb L$-integer is a product of elements of $\mbb L$.
\begin{dfn}A morphism $h:Q\to P$ of fs monoids is \emph{Kummer} (resp. $\mbb
L$-Kummer)
if $h$ is injective
and for every $a\in P$, there exists a positive integer (an $\mbb L$-integer) $n$ such that $a^n\in h(Q)$
(note that if $Q\to P$ is Kummer, $\Spec P\to\Spec Q$ is an homeomorphism).

A morphism $f:X\to Y$ of fs log schemes is said to be \emph{Kummer} (resp. \emph{exact})
if for
every geometric point $\bar x$ of $X$, $\overline M_{Y,f(\bar x)}\to \overline
M_{X,\bar x}$ is Kummer (resp. exact).

A morphism of fs log scheme is \emph{Kummer étale} (or két for short) if it
is Kummer and log étale.\end{dfn}
A morphism $f$ is két if and only if étale locally it is deduced by strict
base change and étale localization from a map $\Spec \mbf Z[P]\to\Spec \mbf
Z[Q]$ induced by a Kummer map $Q\to P$ such that $nP\subset Q$ for some $n$
invertible on $X$.

In fact if $f:Y\to X$ is két, $\bar y$ is a geometric point of $Y$, and
$P\to M_X$ is an exact chart of $X$ at $f(\bar y)$, there is an étale
neighborhood $U$ of $\bar x$ and a Zariski open neighborhood $V\subset
f^{-1}(U)$ of $\bar y$ such that $V\to U$ is isomorphic to $U\times_{\Spec
  \mbf Z[P]}\Spec\mbf Z[Q]$ with $P\to Q$ a $\mbb L$-Kummer morphism where
$\mbb L$ is the set of primes invertible on $U$ (\cite[Prop. 3.1.4]{stix}).

Két morphisms are open and quasi-finite.

\smallskip

The category of két fs log schemes over $X$ (any $X$-morphism
between two such fs log schemes is then két) where the covering families
$(T_i\to T)$ of $T$ are the families that are set-theoretical covering
families (being a set-theoretical covering két family is stable under fs base change) is a site. We
will denote by $X_{\ket}$ the corresponding topos.
Any locally constant finite object of $X_{\ket}$ is representable. 
\begin{dfn} A két fs
log scheme over $X$ which represents such a locally constant finite sheaf
is called a \emph{két cover} of $X$. The category of két covers
of $X_{\ket}$ is denoted by $\KCov(X)$. \end{dfn}
A log geometric point is a log scheme $s$ such that $\mring{s}$ is the
spectrum of a
separably closed field $k$ and $M_s$ is saturated and multiplication
by $n$ on $\overline
M_s$ is an isomorphism for every $n$ prime to the characteristic of $k$.

A \emph{log geometric point} of $X$ is a morphism $x:s\to X$ of log schemes where $s$ is a log
geometric point. A \emph{pointed log scheme} $(X,x)$ is a log scheme $X$ endowed with a log geometric point $x$. A két neighborhood $U$ of $x:s\to X$ in $X$ is a morphism $s\to
U$ of $X$-log schemes where $U\to X$ is két. Then if $x$ is a log geometric
point of $X$, the functor $F_x$ from $X_{\ket}$ to $\Ens$ defined by $\mcal F\mapsto
\varinjlim_U \mcal F(U)$ where $U$ runs through the directed category of két
neighborhoods of $x$ is a point of the topos $X_{\ket}$ and any point of
this topos is isomorphic to $F_x$ for some log geometric point and the family of points $(F_x)$ where $x$ runs through log geometric points of $X$
is a conservative system of points.

\begin{dfn} The inverse limit in the category of saturated
log schemes of the két neighborhoods of $x$ is called the \emph{log strict specialization}, and is denoted by $X(x)$.

If $x$ and $y$ are log
geometric points of $x$, a \emph{specialization} of log geometric points $x\to y$
is a morphism $X(x)\to X(y)$ over $X$.\end{dfn}
A specialization $x\to y$ induces a canonical morphism  $F_y\to F_x$ of functors.

 If there is a specialization $x\to
y$ of the underlying topological points, then there is some specialization $x\to y$ of log
geometric points.

If $X$ is connected, for any log geometric point $x$ of $X$, $F_x$ induces
a fundamental functor $\KCov(X)\to\fSet$ of the Galois category
$\KCov(X)$.

\begin{dfn} The \emph{k\'et fundamental group} $\glog(X,x)$ is the profinite group of automorphisms of
the fundamental functor $\KCov(X)\to\fSet$.\end{dfn}

Strict étale surjective morphisms satisfy effective descent for két
covers (\cite[prop. 3.2.19]{stix}).

If $f:S'\to S$ is an exact morphism of fs log schemes such that $\mring{f}$
is proper, surjective and of finite presentation, then $f$ satisfies
effective descent for két covers (\cite[th. 3.2.25]{stix}).

\begin{prop}\label{covloc} Let $X\to S$ be a morphism of fs log schemes such that $\mring S$ is locally noetherian and $\mring X\to\mring S$ is of finite type. Let $\tilde s$ be a geometric point of $\mring S$ and let $S(\tilde s)$ be the strict localization of $S$ at $\tilde s$ endowed with the pullback log structure. Then the functor $F:\injLim_U\KCov(X_U)\to \KCov(X_{S(\tilde s)})$, where $U$ goes through \'etale neighbrhoods of $\tilde s$, is an equivalence of categories.
\end{prop}
\dem
Let $Y_U$ be a két cover of $X_U$ such that, cofinally on $V$, $Y_V$ is connected. Then $Y_U\times_US(\tilde s)$ is connected according to \cite[prop. 8.4.4]{ega4}. This proves that $F$ is fully faithful, or equivalently the outer morphism of fundamental groups of Galois categories is surjective (\cite[prop. V.6.10]{sga}).

Let $Y\to X_{S(\tilde s)}$ be a két cover. Since one knows that $F$ is fully faithful for any $X$ and surjective \'etale morphisms satisfy effective descent for k\'et covers, one only has to prove the essential surjectivity k\'et locally on $X$, so that one may assume that $X$ has a fs chart $X\to \Spec \mbf Z[P]$. Let $P$ be the characteristic of $\tilde s$. Then there is a $(p')$-Kummer morphism of monoids $P\to Q$ such that $Y_Q:=Y\times_{\Spec \mbf Z[P]}\Spec \mbf Z[Q]$ is strict \'etale over $X_{S(\bar s),Q}:=X_{S(\tilde s)}\times_{\Spec \mbf Z[P]}\Spec \mbf Z[Q]$. There exists a neighborhood $U$ and an \'etale cover $Y_{U,Q}$ of $X_{U,Q}$ such that $Y_Q=Y_{U,Q}\times_{X_U}X_{S(\tilde s)}$. Thus $Y_Q$ is in the essential image of $F$. This proves that the outer morphism of fundamental groups corresponding to $F$ is injective (\cite[cor. V.6.8]{sga}), and therefore $F$ is an equivalence.
\findem

Let us now state the main results to compare log fundamental groups of
different log schemes (in particular specialization comparisons). According to \cite[th. 7.6]{ill}, if $X$ is a log regular fs log scheme, $\KCov(X)$ is equivalent to the category of tamely ramified covers of $X_{\tr}$. If $\mbb L$ is a set of primes invertible on $X$, by taking the pro-$\mbb L$ completion, one gets:

\begin{thm} If $X$ is a log regular fs log scheme and all the primes of $\mbb L$ are
invertible on $X$, then $\KCov(X)^{\mbb L}\to\Covalg(X_{\tr})^{\mbb L}$ is
an equivalence of categories.\end{thm}

For example, if $X$ is a regular scheme and $D$ is a normal crossing divisor and $j:U:=X\backslash D\to X$ is the open immersion, then $M_X=O_X\cap j_*O^*_{X\backslash D}$ is a log structure on $X$ for which $X$ is log regular and $X_{\tr}=U:=X\backslash D$ (for example, if $X=\Spec O_K$ where $O_K$ is a complete discretely valued ring and $D$ is the special point of $X$, then $M_X=O_X\cap j_*O^*_{X\backslash D}$ is the usual log structure of $\Spec O_K$). Thus there is an equivalence of categories $\KCov(X)^{\mbb L}\to\Covalg(U)^{\mbb L}$.

\begin{prop}[{\cite[cor. 2.3]{org2}}] \label{orgsp} Let $S$ be a noetherian strictly local scheme
  with closed point $s$ and let $X$ be a connected fs log scheme such that
  $\mring X$ is proper over $S$. Then
\[\KCov(X)\to\KCov(X_s)\]is an equivalence of categories.\end{prop}
Recall that a strictly local scheme is a henselian scheme such that the residue field at the closed point is separably closed.
One can extend proposition \ref{orgsp} to henselian schemes:

\begin{thm}\label{henselspec}
Let $S$ be a noetherian henselian scheme
  with closed point $s$, and let $X$ be a connected fs log scheme such that
  $\mring X$ is proper over $S$. Then
\[\KCov(X)\to\KCov(X_s)\] is an equivalence of categories.
\end{thm}

\dem
First assume $X_s$ to be geometrically connected. Let $x$ be a log geometric point of $X_s$. Then $X$ is also connected and we have to prove that $\glog(X_s,x)\to\glog(X,x)$ is an isomorphism. Let $\overline s$ be a strict localization of $s$ and let $\overline S$ be the strict localization of $S$ at $\overline s$. Let $\bar x$ be a log geometric point of $X_{\overline s}$ above $x$. Let $S_i$ be a pointed Galois cover of $S$, let $G_i$ be its Galois group and let $s_i=s\times_SS_i$. Then we have a diagram with exact lines:
\[\xymatrix{1 \ar[r] & \glog(X_{s_i},\bar x) \ar[r]\ar[d] & \glog(X_s,x) \ar[d]\ar[r] & G_i \ar[r]\ar[d] & 1\\
1 \ar[r] & \glog(X_{S_i},\bar x) \ar[r] & \glog(X,x) \ar[r] & G_i \ar[r] & 1
}
\]
By taking the projective limit when $S_i$ runs through the category of pointed Galois cover of $S$, one gets a diagram with exact lines
\[\xymatrix{1 \ar[r] & \varprojlim_{S_i}\glog(X_{s_i},\bar x) \ar[r]\ar[d] & \glog(X_s,x) \ar[d]\ar[r] & \ga(S,\overline s) \ar[r]\ar[d] & 1\\
1 \ar[r] & \varprojlim_{S_i}\glog(X_{S_i},\bar x) \ar[r] & \glog(X,x) \ar[r] & \ga(S,\overline s) \ar[r] & 1
}
\]
But, according to proposition \ref{covloc}, $\glog(X_{\overline S},\bar x)\to \varprojlim_{S_i}\glog(X_{S_i},\bar x)$ is an isomorphism.
Similarly $\glog(X_{\overline s},\overline x)\to \varprojlim_{S_i}\glog(X_{s_i},\overline x)$ is an isomorphism. Thanks to proposition \ref{orgsp}, $\glog(X_{\bar s},\overline x)\to\glog(X_{\overline S},\overline x)$ is an isomorphism. Thus $\glog(X_s,x)\to\glog(X,x)$ is also an isomorphism.

In the general case, let $X\to S'$ be the Stein factorization of $X\to S$. For every connected component $S'_j$ of $S'$, let $X_j=X\times_{S'}S'_j$. Since $S'_j$ is henselian and $X_j\to S'_j$ has geometrically connected fibers, one gets that $\KCov(X_{j,s})\to\KCov(X_j)$ is an equivalence of category. Since $\KCov(X)=\prod_j\KCov(X_j)$ and $\KCov(X_s)=\prod\KCov(X_{j,s})$, one gets that $\KCov(X)\to\KCov(X_s)$ is an equivalence of categories.
\findem
\begin{cor}\label{speclogdvr} Let $O_K$ be a complete discretely valued ring endowed with its natural log
structure and let $\mbb L$ a set of prime numbers invertible in $O_K$. Let  $X\to \Spec O_K$ be a proper and log smooth morphism and
let $U:=X_{\tr}\subset X_{\eta}$. There is a natural
equivalence of categories \[\KCov(X_s)^{\mbb L}\simeq\Covalg(U)^{\mbb
  L}.\]\end{cor}
In particular, if $\mring X\to\Spec O_K$ is a semistable model of
$X_{\eta}$, and the log structure on $X$ is given by $M_X=O_X\cap j_*O^*_{X_{\eta}}$ where $j:X_{\eta}\to X$, then $X\to\Spec O_K$ is log smooth and
$X_{\tr}=X_\eta$. We get an equivalence of categories  $\KCov(X_s)^{\mbb L}\simeq\Covalg(X_{\eta})^{\mbb
  L}$, and thus an isomorphism
\[\ga(X_{\eta})^{\mbb L}\to\glog(X_s)^{\mbb L}.\]

Here we recall basic results about saturated morphisms of fs log
schemes. The main reference on the subject is~\cite{sattsuji}, which is
unfortunately unpublished.

\begin{dfn}A morphism of fs monoids $P\to Q$ is \emph{integral} if, for any morphism
of integral monoids $P\to Q'$, the amalgamated sum $Q\oplus_PQ'$ is still
integral.

An integral morphism of fs monoids $P\to Q$ is \emph{saturated} if, for any morphism
of fs monoids $P\to Q'$, the amalgamated sum $Q\oplus_PQ'$ is still
fs.

A morphism $f:Y\to X$ of fs log schemes is \emph{saturated} if for any geometric
point $\bar y$ of $Y$, $\bar M_{X,f(\bar y)}\to\bar M_{Y,\bar y}$ is
saturated.\end{dfn}

If $Y\to X$ is saturated and $Z\to X$ is a morphism of fs log schemes, then
the underlying scheme of $Z\times_XY$ is $\mring Z\times_{\mring X}\mring
Y$.

If $P\to Q$ is a local and integral (resp. saturated) morphism of fs
monoids and $P$ is sharp, the morphism $\Spec \mbf
Z[Q]\to\Spec\mbf Z[P]$ is flat (resp. separable, \emph{i.e.} flat
with geometrically
reduced fibers, cf.~\cite[cor. I.4.3.16]{ogus} and \cite[rem. 6.3.3]{satmorph}).

Let $f:X\to Y$ be log smooth, let $\bar x$ be a geometric point of
$X$ and let $\bar y$ be its image in $Y$. Étale locally on $Y$, there is a chart 
$Y\to\Spec P$ such that $P\to\overline M_{Y,\bar y}$ is an isomorphism. Then, according to~\cite[th. 3.5]{kato},
there is étale locally at $x$ a fs chart $\phi:P\to Q$ of $X\to Y$ such that $Y\to
\Spec \mbf Z[Q]\times_{\mbf Z[P]}X$ is étale such that $\phi$ is injective and the torsion part of $\Coker(\phi^{\gp})$ has order invertible on $X$. Up to localizing $Q$ by the face corresponding to $\bar x$, one can assume that $Q\to M_{X,\bar x}$ is local (and thus exact according to \cite[def. II.2.2.8]{ogus}). Thus if $f$ is integral (resp. saturated), $P\to Q$ is a local and
integral (resp. saturated) morphism of fs monoids and $P$ is sharp. Thus
$f$ is flat (resp. separable).

If $P\to Q$ is an integral morphism of fs monoids, there exists an integer
$n$ such that the pullback $P_n\to Q'$ of $P\to Q$ along
$P\stackrel{n}{\to}P=P_n$ is saturated (theorem~\cite[A.4.2]{satmorph}).

Moreover if $P\to Q$ factors through $Q_0$ such that $P\to Q_0$ is saturated
and $Q_0\to Q$ is $\mbb L$-Kummer, $n$ can be chosen to be an $\mbb
L$-integer.

A morphism $X\to S$ of fs log schemes is said to be \emph{log geometrically saturated} if there exists a két covering $U\to S$ such that $X\times_SU$ is saturated.

For example, if $Y\to S$ is a morphism of fs log schemes, with $\mring S$ locally noetherian and $\mring Y\to \mring S$ of finite type, which factors through $X$ such that $X\to S$ is saturated and $Y\to X$ is két, then $Y\to S$ is log geometrically saturated.

\subsection{Log geometric fundamental groups}

Let $X\to S$ be a morphism of fs log schemes. Let
$\bar x$ be a log geometric point of $X$ and let $\bar s$ be its image in
$S$. The morphism $X\to S$ is said to be \emph{log geometrically connected} at $\bar s$ if there exists a cofinal family of két neighborhoods $U$ of $\bar s$ in $S$ such that $X_U$ is connected. 

The \emph{log geometric fundamental group} of $X$ at $\bar x$ to be
\[\ggeom(X/(S,\bar s),\bar x):=\varprojlim_U \glog(X_U,\bar x),\]
where $U$ runs through k\'et neighborhoods of $\bar s$ in $S$.
If $X\to S$ is log geometrically connected, the category $\ggeom(X/(S,\bar s),\bar x)\tfEns$ of finite sets endowed with a continuous action of $\ggeom(X/s,\bar x)$ is equivalent to the
category \[\KCovgeom(X/(S,\bar s)):=\injLim_{U} \KCov(X_U).\]
In particular, $\ggeom(X/(S,\bar s),\bar x)$ does not depend on $\bar x$ up to outer isomorphism. Therefore, when we work in the category of groups with outer morphisms, the log geometric fundamental group will simply be denoted by $\ggeom(X/(S,\bar s))$.

If $\bar s'\to\bar s$ is a specialization of log geometric points of $S$, there is a natural morphism of pro-log schemes $\text{``}\varprojlim_{\bar s'\in U}\text{''}\ U \to \text{``}\varprojlim_{\bar s\in V}\text{''}\ V$, where $U$ goes through két neighborhoods of $\bar s'$ and $V$ goes through két neighborhoods of $\bar s$. This induces a
functor, 2-functorially in $\bar s'\to \bar s$,
\begin{equation}\label{loggeomspeccompatibility}\KCovgeom(X/(S,\bar s))\to\KCovgeom(X/(S,\bar s')),\end{equation} hence an outer morphism, functorially in $\bar s'\to\bar s$,
\[\ggeom(X/(S,\bar s'))\to\ggeom(X/(S,\bar s)).\]

Let $(S',\bar s')\to (S,\bar s)$ be a morphism of pointed fs log schemes. There is a natural morphism of pro-log schemes $\text{``}\varprojlim_{\bar s'\in U}\text{''}\ U \to \text{``}\varprojlim_{\bar s\in V}\text{''}\ V$, where $U$ goes through két neighborhoods of $\bar s'$ in $S'$ and $V$ goes through két neighborhoods of $\bar s$ in $S$. This induces a functor \[\KCovgeom(X/(S,\bar s))\to\KCovgeom(X'/(S',\bar s'))\]
where $X':=X\times SS'$, hence an outer morphism
\[\ggeom(X'/(S',\bar s'))\to\ggeom(X/(S,\bar s)).\]

\begin{prop}\label{logcovloc}
Let $X\to S$ be a morphism of fs log schemes such that $\mring S$ is locally noetherian and $\mring X\to\mring S$ is of finite type.
Let $\tilde s$ be the geometric point of $\mring S$ defined by $\bar s$ and let $S(\tilde s)$ be the strict localization of $S$ at $\tilde s$ endowed with the pullback log structure. The morphism
\[\ggeom(X_{S(\tilde s)}/(S(\tilde s),\bar s),\bar x)\to \ggeom(X/(S,\bar s),\bar x)\]
is an isomorphism.
\end{prop}
\dem
Let $\mbb L$ be the set of primes invertible at $\bar s$.
By replacing $S$ by an étale neighborhood of $\tilde s$, one can assume that $S$ has a chart $S\to\Spec\mbf Z[P]$ such that the induced map $P\to \overline M_{S,\tilde s}$ is an isomorphism. Extend the map $P^{\gp}\to M^{\gp}_{\bar s,\tilde s}$ into a map $P^{\gp}\otimes\mbf Z[\frac{1}{\mbb L}]\to M^{\gp}_{\bar s,\tilde s}$: this defines for every $\mbb L$-két morphism $P\to Q$ of sharp fs monoids a morphism $\bar s$-log point of
$S(\tilde s)_Q:=S(\tilde s)\times_{\Spec\mbf Z[P]}\Spec\mbf Z[Q]$. When $P\to Q$ goes through $\mbb L$-két morphism $P\to Q$ of sharp monoids, the family $(S(\tilde  s)_Q)$ goes through neighborhoods of $\bar s$ in $S(\tilde s)$.
One has also $\ggeom(X_{S(\tilde s)}/(S(\tilde s),\bar s),\bar x)=\varprojlim_{P\to Q}\varprojlim_{U}\glog(X_U,\bar x)$ where $P\to Q$ goes through $\mbb L$-két morphisms of sharp monoids and $U$ goes through étale neighborhoods of $\bar s$ in $X_Q:=X\times_{\Spec \mbf Z[P]}\Spec\mbf Z[Q]$. Then $\glog(X_{S(\tilde s)_Q},\bar x)\to\varprojlim_U\glog(X_U,\bar x)$, where $U$ goes through étale neighborhoods of $\bar x$ in $X_Q$, is an isomorphism according to proposition \ref{covloc}. Since $\ggeom(X/(S,\bar s),\bar x)=\varprojlim_Q\glog(X_{S(\tilde s)_Q},\bar x)$ one gets the result.
\findem

Assume $\mring S$ to be a henselian local scheme. Let $(T,\bar t)$ be a pointed Galois k\'et cover of $(S,\bar s)$. Then one has an exact sequence:
\[1\to \glog(X_t,\bar x_t)\to\glog(X,\bar x)\to\Gal(t/s),\]
and the right map is onto if $X_t$ is connected.
By taking the projective limit of the previous exact sequence when $(t,\bar
t)$ runs through the directed category of pointed Galois connected
covers of $(s,\bar s)$, one gets an exact sequence
\[1\to \ggeom(X/(S,\bar s),\bar x)\to\glog(X,\bar x)\to\glog(S,\bar s),\]
and the right map is onto if $X\to S$ is log geometrically connected.

Let $O_K$ be a complete discretely valued ring endowed with its natural log
structure and let $\mbb L$ a set of prime numbers invertible in $O_K$. Let  $X\to \Spec O_K$ be a proper and log smooth morphism and
let $U:=X_{\tr}\subset X_{\eta}$. There
is a geometric analog to the specialization isomorphism $\ga(X_{\eta})^{\mbb L}\to\glog(X_s)^{\mbb L}$ of corollary \ref{speclogdvr}:
\begin{thm}[{\cite[th. 1.4]{kisin}}]\label{speclogdvrgeom} There is a natural equivalence of categories
  \[\KCovgeom(X/s)^{\mbb L}\simeq\Covalg(U_{\bar\eta})^{\mbb L}.\]
\end{thm}
It can be deduced from corollary \ref{speclogdvr} thanks to the fact that
 any algebraic cover of $U_{\bar
  \eta}$ is already defined over a tamely ramified extension of $K$ (\cite[prop. 1.15]{kisin}).

\subsection{Specialization of log fundamental groups}\label{speclog}

Let us study specialization of log geometric fundamental groups (that is
the projective limit of the log fundamental groups after taking két
extensions of the base log point).

The only result we will need later on is the following: 
\begin{prop}[{cor. \ref{logsp2}}]\label{logsp}
Let $X\to S$ be a proper and saturated morphism of log schemes such that $\mring S$ is locally noetherian,
  and let $Y\to X$ be a két cover. Let $(s,\bar s)$ and $(s',\bar s')$ be two pointed fs log points of $S$
  and let $\bar s'\to\bar s$ be a specialization of log geometric points. Let $\mbb L$ be a set of primes that does not contain the characteristic of $s$. One has a specialization outer morphism
\[\ggeom(Y_{s'}/(s',\bar s'))^{\mbb L}\to\ggeom(Y_{s}/(s,\bar s))^{\mbb L}.\]

Moreover this morphism factors through
$\ggeom(Y/(S,\bar s))^{\mbb L}$.\end{prop}

To prove this, our main result will
be the invariance of the log geometric fundamental group of an fs
log scheme $X$ saturated and of finite type over an fs log point $S$ with separably closed field by fs
base change that is an isomorphism on the underlying scheme. The
assumptions implies that our base change induces an isomorphism of the
underlying schemes. Working étale locally on this scheme, we are reduced to
the case where this scheme is strictly local, where the log
geometric fundamental group can be explicitly described in terms of the
morphism of monoids $\overline M_{X}\to\overline M_S$.

Combining this base change invariance result with strict base change
invariance of the $\mbb L$-log geometric fundamental group and strict
specialization of the $\mbb L$-log geometric fundamental group (\cite{org2}),
we will get that if $X\to S$ is a proper log smooth saturated morphism, and
$s_2,s_1$ are fs points of $S$ and $\bar s_2\to\bar s_1$ is a
specialization of log geometric points of $S$ over $s_2$ and $s_1$, then
there is a specialization morphism
$\ggeom(X_{s_2})^{\mbb L}\to\ggeom(X_{s_1})^{\mbb L}$.

\begin{lem}\label{lemchgmtbase} Let $s'\to s$ be a strict morphism of fs log points such that
  $\mring s'$ and $\mring s$ are geometric points. Let $\mbb L$ be a set of primes that does not contain the characteristic of $s$. Let $X\to s$ be a
  morphism of fs log schemes such that $\mring X\to\mring s$ is of finite type.

Then $F:\KCov(X)^{\mbb L}\to \KCov(X_{s'})^{\mbb L}$ is an equivalence of categories.\end{lem}
\dem
If $T$ is a connected két cover of $X$, $\mring{T\times_ss'}\to\mring
T\times_{\mring s}\mring s'$ is an isomorphism since $s'\to s$ is
strict. The scheme $\mring T\times_{\mring s}\mring s'$ is connected too, so we get
that the functor $F$ is fully faithful.

As one already knows that $F$ is fully faithful for any $X$, and as strict
étale surjective morphisms satisfy effective descent for két covers, one may prove
the essential surjectivity étale locally, and thus assume that $X$ has a global chart
$X\to\Spec \mbf Z[P]$.

Let $Y'$ be a $\mbb L$-két cover of $X_{s'}$. Then there exists a $\mbb L$-Kummer
morphism of monoids $P\to Q$ such that \[Y'_Q:=Y'\times_{\Spec \mbf
  Z[P]}\Spec\mbf Z[Q]\to X_{s',Q}:=X_{s'}\times_{\Spec \mbf
  Z[P]}\Spec\mbf Z[Q]\] is strict étale (and surjective).

But, since $\mring X_{s',Q}\to\mring X_Q\times_{\mring s}\mring s'$ is an
isomorphism of schemes, $\Covalg(\mring X_{s',Q})^{\mbb L}\to\Covalg(\mring X_Q)^{\mbb L}$ is an equivalence of
categories~(\cite[cor 4.5]{org}). Thus, there is a strict étale cover
$Y_Q$ of $X_Q$ (and thus $Y_Q\to X$ is a két cover) such that $Y'_Q$ is
$X_{s',Q}$-isomorphic to $Y_Q\times_ss'$.

Thus $F$ is an equivalence of categories.
\findem

Let now $s'\to s$ be a morphism of fs log points, such that the underlying
morphism of schemes $\mathring s'\to
\mathring s$ is an isomorphism of geometric points, and let $X\to s$ be a
saturated morphism of fs log schemes with $\mathring X$
noetherian and $\mathring{X}\to\mring{s}$ connected. Since
$X\to s$ is saturated, it is log geometrically connected.

Let $\bar x'$ be a log geometric point of $X'=X\times_ss'$ and let $\bar
x$, $\bar s'$ and $\bar s$ be its image in $X$, $s'$ and $s$ respectively.

\begin{thm}\label{chgmtbase1} Let $s'\to s$ be a morphism of fs log points, such that the underlying
morphism of schemes is an isomorphism of geometric points, and let $X\to s$ be a
saturated morphism of connected noetherian fs log schemes.
The map $\ggeom(X'/(s',\bar s'),\bar x')\to\ggeom(X/(s,\bar s),\bar x)$ is an isomorphism.\end{thm}
\dem
Let $(s_i,\bar s_i)_{i\in I}$ be a cofinal system of pointed Galois
connected két covers of $(s,\bar s)$. Let $\tilde s_i$ be the
reduced subscheme of $s_i$ endowed with the inverse image log structure. Let us write $(X_i,\bar x_i)=(X\times_s\tilde s_i,\bar x\times_{\bar
  s}\bar s_i)$.

Let $(s'_j,\bar s'_j)_{j\in J}$ be a cofinal system of pointed Galois
connected két covers of $(s',\bar s')$. Let $\tilde s'_j$ be the
reduced subscheme of $s'_j$ endowed with the inverse image log
structure. Let us write $(X'_j,\bar x'_j)=(X\times_{s'}s'_j,\bar
x'\times_{\bar s'}\bar s'_j)$.

One has to prove that
 \[\varprojlim_j\glog(X'_j,\bar
 x'_j)\to\varprojlim_i\glog(X_i,\bar
x_i)\]
is an isomorphism, or equivalently that
\[\injLim_i\KCov(X_i)\to\injLim_j\KCov(X'_j)\] is an equivalence of categories.

Since strict \'etale surjective morphisms satisfy effective descent for k\'et covers, the injective limits are filtering and $X$ is quasicompact, it is enough to prove that
\[\injLim_i\KCov(X_i)\to\injLim_j\KCov(X'_j)\] is an isomorphism locally on the \'etale topology of $X$.

According to proposition \ref{covloc}, if $\bar x$ is a geometric point of $\mring X$, then 
\[\injLim_{\bar x\in U}\KCov(U)\to \KCov(X(\bar x)),\]
where $U$ goes through \'etale neighborhoods of $\bar x$, is an equivalence of categories. 
Since $X_i$, $X'_j$ and $X$ have equivalent \'etale topoi, $\bar x$ also defines a point of the \'etale topoi of $X_i$ and $X'_j$. According to~\cite[cor. III.2.1.5.8]{giraud}, one only has to prove that 
\[\injLim_i\KCov(X_i(\bar x))\to\injLim\KCov(X'_j(\bar x))\] for every geometric point $\bar x$ of $X$.

We are thus reduced to the case where $\mathring X$ is a strictly local and
noetherian scheme.
But then (\cite[prop. 3.1.11]{stix}), for $\mathring X$ a strictly local and
noetherian scheme, 
\[\begin{array}{rcl}\varprojlim_i\glog(X_{i,\ket},\bar x_i) & = & \varprojlim_i \Hom(\overline M_{X_i,x_i}^{\gp}, \widehat{\mbf
  Z}^{(p')}) \\ & = & \Hom(\varinjlim_i \overline M_{X_i,x_i}^{\gp},\widehat{\mbf Z}^{(p')})
\\ & = & \Hom(\Coker(\overline M_s^{\gp}\to \overline M_{X,x}^{\gp}),\widehat{\mbf Z}^{(p')}),\end{array}\]
and one has a similar result for $X'$.

Since $\overline M_{X',x'}=\overline M_{X,x}\oplus_{\overline M_s}\overline M_{s'}$, one has $\overline M_{X',x'}^{\gp}=\overline M_{X,x}^{\gp}\oplus_{\overline M_s^{\gp}}\overline M_{s'}^{\gp}$. Thus, $\Coker(\overline M_s^{\gp}\to \overline M_{X,x}^{\gp})\to\Coker(\overline
M_{s'}^{\gp}\to \overline M_{X',x'}^{\gp})$ is an isomorphism.

One thus gets the wanted result.
\findem

Assume now that $(s',\bar s')\to (s,\bar s)$ is a morphism of
pointed fs log points, and that $X\to s$ is log geometrically saturated. Recall that this assumption is satisfied if $X\to s$ goes through $X_0$ such that $X\to X_0$ is két and $X_0\to s$ is saturated.

Let $\mbb L$ be a set of prime that does not contain the characteristic of $s$.
\begin{cor}\label{loginvar}The map of profinite
  groups \[\ggeom(X/(s,\bar s),\bar x)^{\mbb L}\to\ggeom(X'/(s',\bar s'),\bar x')^{\mbb L}\] is an
  isomorphism.\end{cor}
\dem
By replacing $s$ (resp. $s'$) by the closed reduced subscheme of a connected
két cover of $s$ (resp. $s'$), one can assume that $X\to s$ is saturated
($\mring X\to\mring s$ will still be of finite type).

If $(t,\bar t)\to (s,\bar s)$ is a strict étale cover, then $\ggeom(X_t/t,\bar
x_t)\to\ggeom(X/s,\bar x)$ is an isomorphism. Thus, by writing $s_0$ for the
separable closure of $s$ and by taking the
projective limit over pointed strict étale covers (since
$\glog(X_{s_0})=\varprojlim\glog(X_t)$, where $t$ runs through pointed
strict étale covers of $s$), one gets that $\ggeom(X_{s_0}/s_0,\bar
x_0)\to\ggeom(X/s,\bar x)$ is an isomorphism.
One thus may assume that $\mring s$ and $\mring s'$ are geometric points.

Let us consider the fs log scheme $s''$ whose underlying scheme is $\mring
s'$ and whose log structure is the inverse image of the log structure of
$s$.
Thus, one has morphisms $s'\to s''\to s$, where $s'\to s''$ is an
isomorphism on the underlying schemes and $s''\to s$ is strict.
But according to lemma~\ref{lemchgmtbase}, $\glog(X_{s''})^{\mbb L}\to
\glog(X)^{\mbb L}$ and $\glog(s'')^{\mbb L}\to\glog(s)^{\mbb L}$ are isomorphisms. Thus,
\[\ggeom(X_{s''}/s'')^{\mbb L}\to\ggeom(X/s)^{\mbb L}\] is an isomorphism.
By~\ref{chgmtbase1}, $\ggeom(X_{s'}/s')^{\mbb L}\to\ggeom(X_{s''}/s'')^{\mbb L}$ is also an
isomorphism.
\findem

\begin{cor}\label{logsp2}Let $X\to S$ be a proper log geometrically saturated morphism of fs log schemes such that $\mring S$ is locally noetherian. Let $(s,\bar s)$ and $(s',\bar s')$ be two pointed fs log points and
let $\bar s'\to\bar s$ be a specialization of log geometric points. Let $\mbb L$ be a set of primes that does not contain the characteristic of $s$. The functor
\[\phi_{s}:\KCovgeom(X/(S,\bar s))^{\mbb L}\to\KCovgeom(X_s/(s,\bar s))\] is an equivalence.
 Therefore, there is a pair $(\psi_{s/s'},\alpha)$, where $\psi_{s/s'}$ is an exact functor
\[\psi_{s/s'}:\KCovgeom(X_{s}/(s,\bar s))^{\mbb L}\to\KCovgeom(X_{s'}/(s',\bar s'))^{\mbb L}\]
and a natural 2-isomorphism $\alpha$
\[\xymatrix{\KCovgeom(X/(S,\bar s))^{\mbb L} \ar[r]^{\phi_{s/s'}} \ar[d]^{\phi_s} & \KCovgeom(X/(S,\bar s'))^{\mbb L} \ar[d]^{\phi_{s'}} \\ 
  \KCovgeom(X_s/(s,\bar s))^{\mbb L} \ar[r]^{\psi_{s/s'}}\ar@{=>}[ru]^\alpha & \KCovgeom(X_{s'}/(s',\bar s'))^{\mbb L}, } 
\] unique in the sense that if $(\psi_{s/s'}',\alpha')$ satisfies the same conditions, there is a unique 2-isomorphism $\beta:\psi_{s/s'}\to\psi_{s/s'}'$ such that 
$\alpha'\cdot (\phi_s\circ\beta)=\alpha$.
Moreover, if $(s'',\bar s'')$ is a pointed fs log point and $\bar s''\to\bar s'$ is a specialization, then there is a unique isomorphism of functors $\psi_{s/s''}\simeq\psi_{s'/s''}\psi_{s/s'}$ such that the following diagram is 2-commutative:
\[\xymatrix@!0@R=8mm@C=4cm{\KCovgeom(X/(S,\bar s))^{\mbb L} \ar[r] \ar[rrd] \ar[dd] & \KCovgeom(X/(S,\bar s'))^{\mbb L} \ar[dd] \ar[rd] & \\ 
& & \KCovgeom(X/(S,\bar s''))^{\mbb L} \ar[dd]\\
  \KCovgeom(X_s/(s,\bar s))^{\mbb L} \ar[rrd] \ar[r] & \KCovgeom(X_{s'}/(s',\bar s'))^{\mbb L} \ar[rd] &\\
& & \KCovgeom(X_{s''}/(s'',\bar s''))^{\mbb L}
}
\]

\end{cor}
\dem
Let $Z$ be the strictly local scheme of $S$ at $s$ endowed with the inverse
image log structure, and let $z$ be its closed point, endowed with the
inverse image log structure.
The three morphisms
\[\ggeom(X_{s}/(s,\bar s))^{\mbb L}\to\ggeom(X_z/(z,\bar s))^{\mbb L},\]
\[\ggeom(X_z/(z,\bar s))^{\mbb L}\to\ggeom(X_Z/(Z,\bar s))^{\mbb L},\]
\[\ggeom(X_Z/(Z,\bar s))^{\mbb L}\to\ggeom(X/(S,\bar s))^{\mbb L}\]
are isomorphisms according to
corollary~\ref{loginvar}, theorem~\ref{orgsp} and corollary~\ref{logcovloc}. Therefore $\phi_s$ is an equivalence. 
The functor $\psi_{s/s'}$ is then the composition of $\phi_{s'}\phi_{s/s'}$ with a quasi-inverse of $\phi_s$ and the uniqueness is obvious.
The exactness of $\psi_{s/s'}$ comes from the exactness of $\phi_{s/s'}$.
The compatibility with composition is a direct consequence of the uniqueness of $\psi_{s/s''}$.
\findem
Let $\Pt(S)$ be the category whose objects are pointed fs log points $(s,\bar s)$ of $S$, and whose morphisms from $(s,\bar s)$ to $(s',\bar s')$ are specialization of log geometric points $\bar s\to\bar s'$. Corollary \ref{logsp2} tells us that there is a 2-functor $\KCovgeom(X_{(\ )})$ from $\Pt(S)^{\op}$ to the 2-category of Galois categories where 1-morphisms are exact functors which maps $(s,\bar s)$ to $\KCovgeom(X_s/(s,\bar s))$.

This 2-functor induces a functor $\ggeom(X_{(\ )})$ from $\Pt(S)$ to the category of groups with outer morphisms which maps $(s,\bar s)$ to $\ggeom(X_s/(s,\bar s))$.

\section{Cospecialization of graphs of semistable curves}
\subsection{Graphs}
A \emph{graph} $\mbb G$ is given by a set of edges $\mcal E$ a set of vertices
$\mcal V$ and for any $e\in \mcal E$ a set of branches $\mcal B_e$ of
cardinality 2 and a map $\psi_e:\mcal B_e\to\mcal V$. A branch $b$ of $e$ can be
thought of as an orientation of $e$ (or a half-edge), and $\psi_e(b)$ is to be thought of as
the ending of $e$ when $e$ is oriented  according to $b$.

One can also equivalently replace the data of edges and branches of each
edge by the datum of the set of all branches $\mcal B=\coprod_e\mcal B_e$, with an involution
$\iota$ without fixed points (which corresponds heuristically to the reversing of the orientation given by the branch), and a map $\psi:\mcal B\to \mcal V$. The set $\mcal E$ is
then the set of orbits of branches for $\iota$.

A \emph{genuine morphism of graphs} $\phi:\mbb G\to \mbb G'$ is given by a map
$\phi_{\mcal E}:\mcal
E\to \mcal E'$, a map $\phi_{\mcal V}:\mcal V\to\mcal V'$ and for every
$e\in \mcal E$ a bijection $\phi_e:\mcal B_e\to\mcal B'_{\phi_{\mcal
    E}(e)}$ such that the following diagram commutes:
\[\xymatrix{\mcal B_e \ar[d] \ar[r] & \mcal B'_{\phi_{\mcal E}(e)} \ar[d]\\
  \mcal V \ar[r] & \mcal V'}\]
Remark that $\phi_{\mcal E}$ and $\phi_{\mcal V}$ are not enough to define $\phi$ if $\mbb G$ has a loop (\ie an edge whose to branches abut to the same vertex): one can define an automorphism of $\mbb G$ just by inverting the two branches of the loop. Thus, to know how the branches are mapped is important as soon as $\mbb G$ or $\mbb G'$ has loops.

\bigskip

The topological cospecialization for semistable curves will be
given by maps of graphs which are not genuine morphisms.
A \emph{generalized morphism of graphs} $\phi:\mbb G\to\mbb G'$ will be given by:
\begin{itemize}
\item a map $\phi_{\mcal V}:\mcal V\to\mcal V'$;
\item a map $\phi_{\mcal E}:\mcal E\to\mcal E'\coprod\mcal V'$
such that, for any $e\in \mcal E$ such that $\phi_{\mcal E}(e)\in\mcal
  V'$ and for any $b\in\mcal B_e$, $\phi_{\mcal V}\psi(b)=\phi_{\mcal E}(e)$;
\item for any $e\in \mcal E$ such that $\phi_{\mcal E}(e)\in\mcal
  E'$, a bijection $\phi_e:\mcal B_e\to\mcal B'_{\phi_{\mcal
    E}(e)}$ such that the obvious diagram commutes (it is the same diagram
as in the case of genuine morphisms).\end{itemize} 
One can replace the last two data by the data of $\phi_{\mcal B}:\mcal
B\to\mcal B'\coprod\mcal V'$ such that, if $\phi_{\mcal B}(b)\in\mcal B'$,
then $\phi_{\mcal B}(\iota(b))=\iota'(\phi_{\mcal B}(b))$ and $\phi_{\mcal V}\psi(b)=\psi'\phi_{\mcal B}(b)$, and, if $\phi_{\mcal
B}(b)\in\mcal V'$, then $\phi_{\mcal B}(\iota(b))=\phi_{\mcal B}(b)=\phi_{\mcal V}\psi(b)$.

In particular, a genuine morphism is a generalized morphism. Genuine morphisms and generalized morphisms can be composed in an obvious way.

One thus gets a category $\Graph$ of graphs with genuine morphisms and a
category $\GenGraph$ of graphs with generalized morphisms.

\smallskip

There is a geometric realization functor $|\ |:\GenGraph\to\Top$
which maps a graph $\mbb G$ to \[|\mbb G|:=\Coker(\coprod_{b\in\mcal B} \pt_{1,b}\amalg\pt_{2,b}
\rightrightarrows\coprod_{v\in\mcal V} \pt_v\amalg\coprod_{b\in\mcal B}[1/2,1]_b),\]
where \begin{itemize}\item the upper map sends: 
\begin{itemize}\item $\pt_{1,b}$ to $1/2$ in $[1/2,1]_b$ \item
$\pt_{2,b}$ to $1$ in $[1/2,1]_b$,\end{itemize} \item the lower map sends \begin{itemize}\item$\pt_{1,b}$ to
$1/2$ in $[1/2,1]_{\iota(b)}$\item $\pt_{2,b}$ to $\pt_{\psi(b)}$.\end{itemize}\end{itemize}
If $\phi:\mbb G\to\mbb G'$ is a generalized morphism, $|\phi|$ is obtained
by mapping \begin{itemize}\item$\pt_v$ to $\pt_{\phi_{\mcal V}(v)}$,\item $[1/2,1]_b$ to
$[1/2,1]_{\phi_{\mcal B}(b)}$ if $\phi_{\mcal B}(b)\in \mcal B'$ (by the
identity of $[1/2,1]$),\item
$[1/2,1]_b$ to
$\pt_{\phi_{\mcal B}(b)}$ if $\phi_{\mcal B}(b)\in \mcal V'$.\end{itemize}
Remark that, if $\mbb G$ is just a loop, then the geometric realization of the morphism induced by inverting the two branches is not homotopic to the identity: thus $\phi_{\mcal E}$ and $\phi_{\mcal V}$ are not enough in general to characterize the topological behavior of $\phi$.
\subsection{Semistable log curves}
\begin{dfn}\label{stlogcurve}
A morphism $X\to S$ of fs log schemes is a \emph{semistable log curve} if étale
locally on $S$ there is a chart $S\to \Spec P$ such that one of the
following is satisfied:
\begin{itemize}
\item $X\to S$ is a strict smooth curve,
\item $X\to S$ factors through a strict étale morphism $X\to S\times_{\Spec\mbf Z[P]}\Spec\mbf
    Z[Q]$ with $Q=(P\oplus\langle u,v\rangle )/(u + v=p)$ and $p\in P$,
\item $X\to S$ factors through a strict étale morphism $X\to S\times_{\Spec\mbf Z[P]}\Spec\mbf
    Z[P\oplus\mbf N]$.
\end{itemize}
\end{dfn}

A semistable log curve is \emph{strictly semistable} if étale locally on $S$, there
are such maps locally for the Zariski topology of $X$.

\begin{prop}
A morphism $X\to S$ is a semistable log curve if and only if it is a log smooth and
saturated morphism purely of relative dimension 1.
\end{prop}
\dem
The direct sense is obvious. Let $X\to S$ be a saturated log smooth scheme of pure dimension 1. As the definition of a semistable log curve is local for the \'etale topology of $X$ and $S$, one can assume that $S$ has a chart $S\to\Spec\mbf Z[P]$ and $X=S\times_{\Spec\mbf Z[P]}\Spec\mbf Z[Q]$ where $P\to Q$ is an injective local and saturated morphism of monoids, $P$ is sharp and $Q^{\gp}/P^{\gp}$ is invertible on $S$. In particular $T^{\gp}:=\overline Q^{\gp}/\overline P^{\gp}$ is torsionfree. Since $P\to Q$ is saturated, $\Spec\mbf Z[P]\to\Spec\mbf Z[Q]$ is flat and $1=\dim \Spec\mbf Z[P]-\dim \Spec\mbf Z[Q]=\rk P^{\gp}-\rk Q^{\gp}=\rk Q^{\gp}/P^{\gp}\geq \rk T^{\gp}$. Thus $T^{\gp}$ is $\{0\}$ or $\mbf Z$. For every $x\in T$, there exists a unique  $\psi(x)\in \overline Q$ such that $f^{-1}(x)\cup Q=\psi(x)+\overline P$ where $f:\overline Q^{\gp}\to T^{\gp}$ (\cite[prop. I.4.3.14]{ogus}). In particular, if $T^{\gp}=\{0\}$, then $\overline P\to\overline Q$ is bijective, thus $X\to S$ is strict and thus $X\to S$ is smooth.

Assume $T^{\gp }=\mbf Z$. Then $\rk Q^{\gp}/P^{\gp}=\rk T^{\gp}$ and thus $\rk Q^{\gp}=\rk \overline Q^{\gp}$. Since $\overline Q^{\gp}$ is a free abelian group, one can choose a splitting $Q=\overline Q\oplus Q^*$. Since $Q^*\hookrightarrow Q^{\gp}/P^{\gp}$ is finite of order invertible on $S$, $X\to S\times_{\Spec \mbf Z[P]}\Spec\mbf Z[\overline Q]$ is \'etale. Thus one can assume that $Q$ is sharp. Let $T$ be the image of $Q$ in $T^{\gp}$. Then $T=\mbf N$ or $T=\mbf Z$.

First assume $T=\mbf N$. Then $n\psi(1)=\psi(n)+p$ with $p\in P$. Since $P\to Q$ is saturated and $p\leq n\psi(1)$, there exists $p'\in P$ such that $p\leq np'$ and $p'\leq \psi(1)$. Thus, by definition of $\psi(1)$, $p'=0$, thus $p=0$ and $\psi(n)=n\psi(1)$. Thus $Q=P\oplus \mbf N\psi(1)$.

If $T=\mbf Z$, let $u=\psi(1)$ and $v=\psi(-1)$. Since $\psi(u+v)=0$, $p:=u+v\in P$. As in the previous case, if $n\geq 0$, then $\psi(n)=n\psi(1)$ and $\psi(-n)=n\psi(-1)$. Thus $Q=P\oplus\langle u,v\rangle /(u+v=p)$.
\findem
The underlying morphism of schemes $\mring X\to \mring S$ is a semistable
curve. In particular, if $\mring S$ is a geometric point, one can associate to $X$ a graph $\mbb G(X)$ in the
following way: the vertices are the irreducible components of $X$,
the edges are the nodes. If $x$ is a node, then the henselization $X(x)$ of $X$ at $x$ has two irreducible components: these components are the branches of the edge corresponding to $x$. If $z$ is an irreducible component of $X(x)$ and $z'$ is the irreducible component of $X$ containing the image of $z$ in $X$, the branch corresponding to $z$ abuts to the vertex corresponding to $z'$ (this graph does not depend of the
log structure).

\smallskip

If $X\to S$ is a proper semistable log curve and $X'\to X$ is a két
cover, then for any log geometric point $\bar s$ of $S$, there is a két
neighborhood $U$ of $\bar s$ such that $X'_U\to U$ is saturated. Then
$X'_U\to U$ is also a semistable curve.

The morphism $\mring X'_{\bar s}\to\mring X_{\bar s}$ induces a
genuine morphism $\mbb G(X'_{\bar s})\to\mbb G(X_{\bar s})$ of graphs.

\subsection{Topological cospecialization of semistable curves}

Let $f:X\to S$ be a semistable curve such that $S$ is locally noetherian, and let $\bar s_2\to \bar s_1$ be a specialization of
geometric points of $S$. In this section we will define a cospecialization
map of graphs $\mbb G(X_{\bar s_1})\to \mbb G(X_{\bar s_2})$.

\begin{prop}\label{lemcospstr}Let $f:X\to S$ be a strictly semistable curve such that $S$ is strictly local and noetherian. Let $s_1$ be the closed point of $S$, and let $s_2$ be a point of $S$. Let $x$
  be a node or a generic point of $X_{s_1}$. Let $X(x)$ be the localization of $X$ at $x$. Then $X(x)_{s_2}$ is either
  contained in the smooth locus of a geometrically irreducible component, denoted by $F(x)$, of
  $X(x)_{s_2}$ or contains a single node, denoted by $F(x)$, of $X(x)_{s_2}$, which is rational.
\end{prop}
\dem

Let $A$ be the noetherian strictly local ring such that $S=\Spec A$. By replacing $S$ by the closure of $s_2$ endowed with the reduced scheme structure, one can assume that $s_2$ is the
generic point of $S$ and $S$ is integral. Indeed, nonempty closed subschemes of henselian schemes are henselian (\cite[cor. to \S{} 3.prop 2]{hensel}) and keep the same residue field at the special point, therefore the closure of $s_2$ is a strictly local scheme.

\begin{enumerate}[(i)]
\item
If $x$ is in the smooth locus of $X_{s_1}$, $X\to S$ is smooth at $x$, and
$X(x)_{s_2}$ is geometrically connected by local 0-acyclicity of smooth
morphisms.
\item
If $x$ is a node, one can assume that $f$ factors
through an étale morphism $X\to\Spec B$ with $B=A[u,v]/(uv-a)$ and
$a(s_1)=0$.

If $a=0$ let $Z=X\times_{\Spec B}\Spec A$ where $g:B\to A$ is defined by
$g(u)=g(v)=0$ (this is the closed subscheme of $X$ defined by the
node; in particular $Z_{s_2}$ is the union of all the nodes of $X_{s_2}$). The morphism $Z\to S$ is étale and thus $Z(x)\to S$ is an isomorphism. Thus
$Z(x)_{s_2}$ is just a rational point $F(x)$.

If $a\neq 0$, then $a(s_2)\neq 0$ and thus $X_{s_2}$  is smooth. Since $X\to S$
is a semistable curve, it is separable (\ie flat with separable geometric
fibers). Let $B:=O_{X,x}$ be the noetherian local ring such that $X(x)=\Spec B$. By applying \cite[cor. 18.9.8]{ega4} to $\Spec B=X(x)\to \Spec A=S$, one gets that $X(x)_{s_2}$ is
geometrically connected.
\end{enumerate}
\findem

Let $f:X\to S$ be a semistable curve, and let $\bar s_2\to \bar s_1$ be a specialization of
geometric points of $S$. One can apply lemma \ref{lemcospstr} to $X_{S(\bar s_1)}\to S(\bar s_1)$ and to the Zariski point $s_2$ corresponding to $\bar s_2$. Let $x$ be a node or a generic point of $X_{\bar s_1}$.
 If $F(x)$ is a rational node of $X_{s_2}$, then it defines an edge $F_0(x)$ of $\mbb G_{\bar s_2}$. If $F(x)$ is a geometrically irreducible component of $X_{s_2}$, then it defines a vertex $F_0(x)$ of $\mbb G_{\bar s_2}$.

\begin{lem}\label{chgbase}Let $S'=\Spec A'\to S=\Spec A$ be a local morphism of noetherian strictly local schemes. Let $s'_1$ be the closed point of $S'$ and let $s'_2$ be a point of $S'$ above $s_2$. Let $X'=X\times_SS'$. Let $x'\in X'_{s'_1}$ be above $x$. Let $F'$ be defined analogously to $F$ but with the curve $X'\to S'$. Then $F(x)$ is the image of $F'(x')$ by the map $X'\to X$. 
\end{lem}
\dem
Let $z$ be the image of $F'(x')$ by the map $X'\to X$.
Since $x'$ is in the closure of $F'(x')$, $x$ is in the closure of $z$, and therefore $F(x)$ is in the closure of $z$.
One only has to prove that if $F(x)$ is a
node, then $z$ is also a node. One can assume that $X\to S$ factors through an \'etale morphism $X\to\Spec B$ with $B=A[u,v]/(uv-a)$ and
$a(s_2)=0$. Then $X'\to S'$ factors through the \'etale map $X\to \Spec B'$ with $B'=A'[u,v]/(uv-a')$ with $a'(s'_2)=0$, and thus $F'(x')$ and $z$ are also nodes.
\findem

\begin{lem}\label{chgopen}
If $\phi:X'\to X$ is a quasifinite open morphism of strictly semistable curves over
$S$ which maps nodes to nodes on every fiber, then $\phi F_0'=F_0\phi$.\end{lem}
For example, the assumption is satisfied if $X'\to X$ is \'etale or if $X'\to X$ is a k\'et morphism of strictly semistable log curves.
\dem 
One can assume that $S$ is strictly local with closed point $\bar s_1$, so that one has to prove that $\phi F=F\phi$.
Since
$\phi(x)$ is in the closure of $\phi F'(x)$, $F\phi(x)$ is in the
closure of $\phi F'(x)$. One only has to prove that if $F\phi(x)$ is a
node, then $\phi F'(x)$ is also a node. Let us assume that  $F \phi(x)$
is a node of $X_{s_2}$. Let $z_1$ and $z_2$ be the two generic points of
the irreducible components of $X_{s_2}$ whose closures contain $F \phi(x)$
(and thus also $\phi(x)$). Since $\phi$ is open, there exists $z'_1$ and
$z'_2$ in $X(x)_{s_2}$ such that $\phi(z'_1)=z_1$ and
$\phi(z'_2)=z_2$. Thus $X(x)_{s_2}$ cannot be in a single irreducible
component of $X_{s_2}$, and thus $F'(x)$ is a node of $X'_{s_2}$. By
assumption, $\phi F'(x)$ is a node of $X_{s_2}$.\findem

\begin{prop}\label{cosptop}
Let $S$ be a locally noetherian scheme and let $\bar s_2\to\bar s_1$ be a specialization of geometric points.
There is a unique way to associate to every semistable curve $X\to S$ a generalized morphism of graphs \[\psi:\mbb G(X_{\bar
  s_1})\to\mbb G(X_{\bar s_2})\]
\begin{itemize}
\item  which is functorial for étale morphisms $X'\to X$,
\item such that if $f:X\to S$ is strictly semistable, $\psi(x)=F_0(x)$ for any
node or generic point $x$ of $\mbb G(X_{\bar s_1})$.
\end{itemize}
\end{prop}
\dem
 After replacing $S$ by its strict localization at $\bar s_1$, one can assume that $S$ is strictly local and $\bar s_1$ is the closed point.

First, let us prove the uniqueness. Let $X\to S$ be a semistable curve. Let $x$ be a node or a vertex of $\mbb G(X_{\bar s_1})$. Let $X'\to X$ be a surjective \'etale morphism such that $X'\to S$ is strictly semistable. Let $x'$ be a preimage of $x$ in $\mbb G(X'_{\bar s_1})$. Then, by functoriality, $\psi(x)$ must be the image of $F_0(x')$ by the map $\mbb G_{X'_{\bar s_2}}\to\mbb G_{X_{\bar s_2}}$. Moreover if $b$ is a branch of $\mbb G(X_{\bar s_1})$, let $b'$ be a preimage in $\mbb G(X'_{\bar s_1})$. Since $\mbb G_{X_{\bar s_2}}$ has no loop, $\psi(b')$ is uniquely defined by the vertex it is ending at, $\psi(b)$ is the image of $\psi(b')$ by the map $\mbb G_{X'_{\bar s_2}}\to\mbb G_{X_{\bar s_2}}$. This proves the uniqueness.

\bigskip

Let us now construct $\psi$.

Let $f:X\to S$ be a strictly semistable curve, and let $\bar s_2\to \bar s_1$ be a specialization of
geometric points of $S$.

If $e$ is a vertex of $\mbb G(X_{\bar s_1})$, then $\psi(e):=F_0(x)$ where $x$ is the node of $X_{\bar s_1}$ corresponding to $e$.
If $v$ is a vertex of $\mbb G(X_{\bar s_1})$, then $\psi(v):=F_0(x)$ where $x$ is the generic point of the irreducible component of $X_{\bar s_1}$ corresponding to $e$. Let $b$ be a branch
  of an edge $e$ in $\mbb G(X_{s_1})$ that abuts to a vertex $v$. Then $F(x)\subset F(s)$, where $x$ is the node corresponding to $e$ and $s$ is the generic point of the irreducible component corresponding to $v$. If $F(x)=F(s)$, then $\psi(b):=F_0(x)=F_0(s)$. Otherwise, $\psi(e)$ is an edge and $\psi(v)$ is a vertex, and there is a branch $b'$ of $\psi(e)$ abutting to $\psi(v)$. Since $X_{\bar s_2}$ is strictly
  semistable, this branch is unique. Let $\psi(b)=b'$.

The compatibility with \'etale morphisms is a direct consequence of lemma \ref{chgopen}.

\bigskip

If $X\to S$ is now a general semistable curve, one choses a surjective étale morphism
 $X'\to X$, and let $X''=X'\times_XX'$.

One has a commutative diagram with genuine lines
\begin{equation}\label{cokergrph}\xymatrix{\mbb G(X''_{\bar s_1}) \ar[d]^{\psi''}\doublear[r] & \mbb
  G(X'_{\bar s_1}) \ar[r]\ar[d]^{\psi'} &  \mbb G(X_{\bar s_1}) \\
\mbb G(X''_{\bar s_2}) \doublear[r] & \mbb
  G(X'_{\bar s_2}) \ar[r] & \mbb G(X_{\bar s_2})}\end{equation}
such that \[\mcal V_{\mbb G(X_{\bar s_2})}=\Coker(\mcal V_{\mbb G(X''_{\bar s_2})} \rightrightarrows  \mcal V_{\mbb
  G(X'_{\bar s_2})})\]
and \[\mcal E_{\mbb G(X_{\bar s_2})}=\Coker(\mcal E_{\mbb G(X''_{\bar s_2})} \rightrightarrows  \mcal E_{\mbb
  G(X'_{\bar s_2})}).\]
By taking the cokernel one gets maps $\psi_{\mcal V}:\mcal V_{\mbb G(X_{\bar s_1})}\to \mcal V_{\mbb G(X_{\bar s_2})}$ and $\psi_{\mcal E}:\mcal E_{\mbb G(X_{\bar s_1})}\to \mcal E_{\mbb G(X_{\bar s_2})}\coprod\mcal V_{\mbb G(X_{\bar s_2})}$.

 Let $e$ be an edge of $\mbb G(X_{\bar s_1})$ such that $\psi_{\mcal E}(e)\in \mcal E_{\mbb G(X_{\bar s_2})}$.
Let $e'$ be an edge of $\mbb G(X_{\bar s_1})$ mapping to $e$. One has bijections $\mcal B_e\leftarrow\mcal B_{e'}\to\mcal B_{\psi'_{\mcal E}(e')}\to\mcal B_{\psi_{\mcal E}(e)}$, hence a bijection $\psi_e:\mcal B_e\to\mcal B_{\psi_{\mcal E}(e)}$. Let $e'_1$ and $e'_2$ be edges of $\mbb G(X_{\bar s_1})$ mapping to $e$. There exists an edge $e''\in \mcal E_{\mbb G(X_{\bar s_2})}$ mapping to $e'_1$ and $e'_2$ by the two maps $\mbb G(X''_{\bar s_1}) \to \mbb
  G(X'_{\bar s_1})$. One gets a commutative diagram of bijections:
\[\xymatrix{ & \mcal B_{e_1'} \ar[dl] \ar[r] & \mcal B_{\psi'_{\mcal E}(e_1')} \ar[dr]& \\
\mcal B_e & \mcal B_{e''} \ar[d]\ar[u]\ar[r]\ar[l] & \mcal B_{\psi''_{\mcal E}(e'')} \ar[d]\ar[u]\ar[r] & \mcal B_{\psi_{\mcal E}(e)}, \\
& \ar[ul] \ar[r] \mcal B_{e_2'} & \mcal B_{\psi'_{\mcal E}(e_2')} \ar[ur] & 
}
\]
which proves that the bijection $\psi_e$ does not depend on the choice of $e'$. The wanted compatibilities between $\psi_{\mcal E}$, $\psi_{\mcal V}$ and $\psi_e$ come directly from the corresponding compatibilities between $\psi'_{\mcal E}$, $\psi'_{\mcal V}$ and $\psi'_{e'}$.
Therefore, there is a unique generalized morphism of graphs $\psi:\mbb G(X_{\bar
  s_1})\to\mbb G(X_{\bar s_2})$ making the diagram \eqref{cokergrph} commutative. 

This morphism $\psi$ does not depend of the choice of $X'$. Indeed let $X'_1\to X$ and $X'_2\to X$ be two surjective \'etale morphisms such that $X'_1$ and $X'_2$ are strictly semistable. By considering $X'_1\times_XX'_2\to X$, one can assume that there is a $X$-morphism $X'_2\to X'_1$. Then one has a diagram
\[\cube{\mbb G(X'_{2,\bar s_1})}{\mbb G(X'_{1,\bar s_1})}{\mbb G(X'_{2,\bar s_2})}{\mbb G(X'_{1,\bar s_2})}{\mbb G(X_{\bar s_1})}{\mbb G(X_{\bar s_1})}{\mbb G(X_{\bar s_2})}{\mbb G(X_{\bar s_2})}\]
where the horizontal maps of the lower square are identities and the frontward maps of the lower square are the two versions of $\psi$ defined in terms of $X'_1$ and $X'_2$. Since the upper face and the vertical faces are commutative and the vertical maps are surjective, the lower square is also commutative. Therefore $\psi$ does not depend on $X'$.

Let us show the functoriality of $\psi$ with respect to \'etale morphisms. Let $X_2\to X_1$ be an étale morphism. Let $X_1'\to X_1$ be a surjective \'etale morphism such that $X'_2\to S$ is strictly semistable. Let $X'_2:=X'_1\times_{X_1}X$. Consider the diagram
\[\cube{\mbb G(X'_{2,\bar s_1})}{\mbb G(X'_{1,\bar s_1})}{\mbb G(X'_{2,\bar s_2})}{\mbb G(X'_{1,\bar s_2})}{\mbb G(X_{2,\bar s_1})}{\mbb G(X_{1,\bar s_1})}{\mbb G(X_{2,\bar s_2})}{\mbb G(X_{1,\bar s_2}).}\]
Since the upper face and the vertical faces are commutative and the vertical maps are surjective, the lower square is commutative.

\findem

\begin{prop}
Let $f:S'\to S$ be a morphism of locally noetherian schemes. Let $\bar s'_2\to\bar s'_1$ be a specialization of geometric points of $S'$, and let $\bar s_2\to\bar s_1$ be the image in $S$. Let $X\to S$ be a semistable curve and let $X'=X\times_{S'}S$. Then the diagram
\begin{equation}\xymatrix{\mbb G_{X'_{\bar s'_1}}\ar@{=}[d]\ar[r]^{\psi'} & \mbb G_{X'_{\bar s'_2}}\ar@{=}[d] \\
\mbb G_{X_{\bar s_1}}\ar[r]^{\psi} & \mbb G_{X_{\bar s_2}},
 }
\end{equation}
where $\psi$ and $\psi'$ are the cospecialization maps, is commutative.
\end{prop}
\dem
Up to replacing $S'$ by its strict localization at $\bar s'_1$ and $S$ by its strict localization at $\bar s_1$, one can assume that $S'\to S$ is a local morphism of strictly local schemes and that $\bar s'_1$ and $\bar s_1$ are the closed points of $S'$ and $S$. 
Let $\psi_0$ be the composition $\mbb G_{X_{\bar s_1}}=\mbb G_{X'_{\bar s'_1}}\stackrel{\psi'}{\to}\mbb G_{X'_{\bar s'_2}}=\mbb G_{X_{\bar s_2}}$. Since $\psi'$ is compatible with \'etale morphisms, $\psi_0$ is also compatible with étale morphisms. Let $f:X\to S$ be a strictly semistable morphism and let $x$ be a node or a vertex of $\mbb G_{X_{\bar s_1}}$. Let $x'$ be the corresponding node or vertex of $\mbb G_{X'_{\bar s'_1}}$. Then $fF_0(x')=F_0(x)$ according to lemma \ref{chgbase}. Therefore $\psi_0(x)=f\psi'(x')=fF_0(x')=F_0(x)$. Therefore, by uniqueness in proposition \ref{cosptop}, one has $\psi_0=\psi$.
\findem

\begin{prop}\label{specializationcompatibility}
Let $X\to S$ be a semistable curve. Let $\bar s_3\to\bar s_2$ and $\bar s_2\to\bar s_1$ be specializations of geometric points of $X$. Then the diagram
\[\xymatrix{\mbb G_{X_{\bar s_1}}\ar[r]^{\psi_{12}}\ar[rd]^{\psi_{13}} & \mbb G_{X_{\bar s_2}} \ar[d]^{\psi_{23}}\\
& \mbb G_{X_{\bar s_3}},
 }
\]
where $\psi_{12}$, $\psi_{13}$ and $\psi_{23}$ are cospecialization maps, is commutative.
\end{prop}
\dem
The morphism $\psi_{23}\psi_{12}:\mbb G_{X_{\bar s_1}}\to\mbb G_{X_{\bar s_3}}$ is functorial with respect to étale morphisms $X'\to X$. By uniqueness in proposition \ref{cosptop}, it is enough to prove that $\psi_{23}\psi_{12}(x)=\psi_{13}(x)$ for every node or edge $x$ of $\mbb G_{X_{\bar s_1}}$ assuming that $X$ is strictly semistable. Since $\psi_{23}\psi_{12}(x)$ specializes to $x$, $\psi_{13}(x)$ is in the closure of $\psi_{23}\psi_{12}(x)$. Therefore, one only has to prove that if $\psi_{13}(x)$ is a node, $\psi_{23}\psi_{12}(x)$ is also a node. Then up to \'etale localization, one can assume $S=\Spec A$ and $X\to \Spec A$ goes through an étale morphism $X\to \Spec B$ where $B=A[u,v]/(uv-a)$ with $a(\bar s_3)=0$, in which case it is obvious.
\findem

\begin{prop}\label{chgopen2}
If $\phi:X'\to X$ is a quasifinite open morphism of semistable curves over
$S$ which maps nodes to nodes on every fiber, then $\phi \psi'=\psi\phi$, where $\psi:\mbb G_{X_{\bar s_1}}\to\mbb G_{X_{\bar s_2}}$ and $\psi':\mbb G_{X'_{\bar s_1}}\to\mbb G_{X'_{\bar s_2}}$ are cospecialization maps.\end{prop}
\dem
Since the lemma is true if $X'\to X$ is étale, one only has to prove it locally on $X'$ and $X$ for the étale topology. Therefore one can assume that $X$ and $X'$ are strictly semistable curves over $S$. According to lemma \ref{chgopen}, for any node or edge $x$ of $\mbb G_{X'_{\bar s_1}}$, $\phi\psi'(x)=\phi F'_0(x)=F_0\phi(x)=\psi\phi(x)$. Since $\mbb G_{X'_{\bar s_1}}$ has no loop, this implies that $\phi\psi'=\psi\phi$.
\findem

\smallskip

We want to know when this generalized morphism of graphs is an isomorphism.

\begin{prop}\label{isomproper}
Keeping the notations of proposition \ref{cosptop}, if $\psi:\mbb G(X_{\bar s_1})\to\mbb G(X_{\bar s_2})$ is a genuine morphism
of graphs and $f$ is proper, then $\psi$ is an isomorphism.
\end{prop}
\dem
One may assume $S=\Spec A$ to be strictly local and integral with special
point $s_1$ and generic point $s_2$.
The assumption means that étale locally on the special fiber (and thus on
$X$ by properness), $X$ is isomorphic to $\Spec A[u,v]/uv$ or is smooth.

Let $Z\subset X$ be the non smooth locus
of $X\to S$, endowed with the reduced subscheme structure. $Z\to S$ is
étale (as can be seen étale locally over $X$), and proper. One thus gets
that $F$ induces a bijection between nodes of $X_{\bar s_1}$ and $X_{\bar
  s_2}$.

Let $\widetilde X$ be the blowup of $X$ along $Z$. 
When $X=\Spec A[u,v]/(uv)$, $Z$ is defined by the ideal generated by $u$
and $v$, and $\widetilde X=\Spec A[u]\coprod \Spec A[v]$.
Thus by looking étale locally over $X$, one sees
that $\widetilde X$ is smooth over $S$, and that $\widetilde X_s$ is simply the
normalization of $X_s$.
Since we assumed $X\to S$ to be proper, $\widetilde X\to S$ is smooth and
proper, thus its Stein factorization induces a bijection between the
connected components of $\widetilde X_{\bar s_1}$ and $\widetilde X_{\bar
  s_2}$, and thus the map between the irreducible components of $
X_{\bar s_1}$ and $X_{\bar s_2}$ is a bijection too.
\findem

\begin{prop}\label{invstrates}
Let $f:X\to S$ be a log semistable curve and let $\bar s_2\to
\bar s_1$ be a specialization of log geometric point.

Assume $\overline M_{S,\bar s_1}\to\overline M_{S,\bar s_2}$ is an isomorphism.
Then $\psi:\mbb G(X_{\bar s_1})\to\mbb G(X_{\bar s_2})$ is a genuine morphism of
graphs. \end{prop}
\dem
One can assume $S$ to be strictly local, integral with generic point $s_2$:
$S=\Spec A$, with a chart $P\to A$.

To show that it is a genuine morphism, one only has to prove that $\psi(e)$
is an edge if $e$ is an edge of $\mbb G(X_{s_1})$. This is not changed by
an étale morphism, so that one can simply assume $X=\Spec A\otimes_{\mbf
  Z[P]}\mbf Z[Q]$ with $Q=(P\oplus\langle u,v\rangle )/(u+v=p)$ and $p\in P$, such that
the image of $p$ in $M_{\bar s_1}$ is not invertible. Thus the image of $p$
in $M_{\bar s_2}$ is not invertible and thus $X=\Spec A[u,v]/(uv=0)$, which
gives the wanted result.
\findem

\subsection{Topological cospecialization and két morphisms}

\begin{prop}\label{topcospket}Let $S$ be a fs log scheme such that $\mring S$ is locally noetherian. Let $(s_2,\bar s_2)$ and $(s_1,\bar s_1)$ be two pointed fs log points of $S$ and let $\bar s_2\to \bar s_1$ be a specialization of log geometric points. Let $X\to S$ be a proper semistable log curve and let $Y_{\bar s_2}$ be a log geometric két cover of $X_{s_1}/(s_1,\bar s_1)$. There is a unique morphism of graphs
\[\phi:\mbb G(Y_{\bar s_1}) \to \mbb G(Y_{\bar
  s_2}),\]
where $Y_{\bar s_2}$ is the image of $Y_{\bar s_1}$ by the functor $\KCov(X_{\bar s_1})\to\KCov(X_{\bar s_2})$ given by corollary~\ref{logsp2}, such that, if $U$ is a két neighborhood of $\bar s_1$ in $S$ and $Z\to X_U:=X\times_SU$ is an extension of $Y_{\bar s_1}$ such that $Z\to U$ is saturated (and therefore a semistable log curve), then the diagram
\[
 \xymatrix{\mbb G_{Z_{\bar s_1}}\ar@{=}[d]\ar[r]^{\psi} & \mbb G_{Z_{\bar s_2}}\ar@{=}[d]\\
\mbb G_{Y_{\bar s_1}}\ar[r]^{\phi} & \mbb G_{Y_{\bar s_2}},
}
\] where $\psi$ is the cospecialization morphism defined by proposition \ref{cosptop}, is commutative.

Moreover $\phi$ is functorial with respect to morphisms $Y'_{\bar s_1}\to Y_{\bar s_1}$ of log geometric két covers of $X_{s_1}/(s_1,\bar s_1)$ and with respect to composition of specializations of log geometric points.

If $\overline M_{S,s_1}\to\overline M_{S,s_2}$ is an isomorphism, then $\phi$ is an isomorphism.

\end{prop}
\dem
According to corollary~\ref{logsp2}, there exists a két neighborhood $U$ of $\bar s_1$ and a két cover $Z\to X_U$ which extends $Y_{\bar s_2}$. Up to replacing $U$ by a smaller két neighborhood, one can assume that $Z\to U$ is saturated. This proves the uniqueness. One only has to prove that the morphism $\phi$ one gets does not depend on the choice of $U$ and $Z\to X_U$. Let $U$ and $U'$ be two k\'et neighborhoods of $\bar s_1$ and let $Z\to X_U$ and $Z'\to X_{U'}$ be két covers that extend $Y_{\bar s_2}$. Since $\KCovgeom(X_{s_1}/(s_1,\bar s_1))\to\KCovgeom(X/(Z,\bar s_1))$ is an equivalnce there exists a k\'et neighborhood $U''$ of $\bar s_1$ in $U\times_SU'$ and an isomorphism $Z'_{U''}\simeq Z_{U''}$. Therefore one can assume that there is a morphism $U'\to U$ over $S$ and that $Z'=Z_{U'}$. Since the specializations maps of proposition~\ref{cosptop} are compatible with a base change $U'\to U$, $U$ and $U'$ define the same morphism $\phi$. This proves the existence of $\phi$.

Let $Y'_{\bar s_1}\to Y_{\bar s_1}$ be a morphism of log geometric két covers of $X_{s_1}/(s_1,\bar s_1)$. There exists a két neighborhood $U$ of $\bar s_1$ and extensions $Z\to X_U$ and $Z'\to Z$ of $Y_{\bar s_1}$ and of $Y'_{\bar s_1}\to Y_{\bar s_1}$ such that $Z'\to U$ is saturated. The compatibility of $\phi$ with $Y'_{\bar s_1}\to Y_{\bar s_1}$ is equivalent to the compatibility of $\psi$ with $Z'\to Z$, which is given by proposition \ref{chgopen2}.

Let $(s_3,\bar s_3)$ be a pointed fs log point and let $\bar s_3\to \bar s_2$ be a specialization. Let $U$ be a két neighborhood of $\bar s_1$ and $Z\to X_U$ be an extension of $Y_{\bar s_1}$ such that $Z\to U$ is saturated. The compatibility of $\phi$ with the composition of specializations for $Y_{\bar s_1}$ is equivalent to the compatibility of $\psi$ with the composition of specialization for $Z$, which is given by lemma \ref{specializationcompatibility}.

If $\overline M_{S,s_1}\to\overline M_{S,s_2}$ is an isomorphism,
$\overline M_{U,s_1}\to\overline M_{U,s_2}$ is still an isomorphism, so
that one can still apply proposition~\ref{invstrates} to an extension $Z\to U$ of $Y_{\bar s_1}$: the morphism $\phi$ is a genuine morphism of graphs. According to proposition~\ref{isomproper}, $\phi$ is an isomorphism.
\findem
If $Y'_{\bar s_1}\to Y_{\bar s_2}$ is a morphism of log geometric két covers, then the following diagram is commutative:
\[\begin{array}{ccc} \mbb G(Y'_{\bar s_1}) & \to  & \mbb G(Y'_{\bar s_2})\\
\dar & & \dar \\ \mbb G(Y_{\bar s_1}) & \to & \mbb G(Y_{\bar
  s_2})\end{array}\]
If $\overline M_{S,s_1}\to\overline M_{S,s_2}$ is an isomorphism,
$\overline M_{U,s_1}\to\overline M_{U,s_2}$ is still an isomorphism, so
that one can still apply proposition~\ref{invstrates} to $Y$: the morphism $\mbb G(Y_{\bar s_1}) \to \mbb G(Y_{\bar
  s_2})$ is a genuine morphism of graphs.

\section{Cospecialization of tempered fundamental groups}
 \subsection{Tempered fundamental groups}
Let $K$ be a complete nonarchimedean field.

Let $\mbb L$ be a set of prime numbers (for example, we will denote by
$(p')$ the set of all primes except the residual characteristic $p$ of
$K$). An $\mbb L$-integer will be an integer which is a product of elements
of $\mbb L$.

If $X$ is a $K$-algebraic variety, $X^{\an}$ will be the
$K$-analytic space in the sense of Berkovich associated to $X$.

A morphism $f:S'\to S$ of analytic spaces is said to be an \emph{étale cover} if $S$ is covered by open subsets $U$ such that $f^{-1}(U)=\coprod V_j$ and $V_j\to U$ is étale finite~(\cite{dJ1}).

For example, étale $\mbb L$-finite covers (\emph{i.e.} finite étale
covers that are dominated by a Galois cover $S''$ of $S$ such that
$\#\Gal(S''/S)$ is an $\mbb L$-integer), also called \emph{$\mbb
  L$-algebraic covers}, and covers in the usual topological sense for
the Berkovich topology, also called \emph{topological covers}, are
étale covers.

Then, André defines tempered covers in \cite[def. 2.1.1]{andre1}. We generalize this definition to
$\mbb L$-tempered covers as follows:
\begin{dfn} \label{def:rvt:temp}
An étale cover $S' \to S$ is \emph{$\mbb L$-tempered} if it is a
 quotient of the composition of a topological cover $T' \to T$ and of a
 $\mbb L$-finite étale cover
 $T \to S$.
\end{dfn}
This is equivalent to say that it becomes a topological cover after
pullback by some $\mbb L$-finite étale cover.

Let $X$ be a $K$-analytic space. We denote by $\CovtempL(X)$ (resp. $\Covalg(X)^{\mbb L}$,
$\Covtop(X)$) the category of $\mbb L$-tempered covers (resp. $\mbb L$-algebraic covers, topological covers) of $X$ (with the obvious morphisms).

\bigskip

A geometric point of a $K$-analytic space $X$ is a morphism of Berkovich spaces $\mcal M(\Omega)\to X$ where $\Omega$ is an algebraically closed complete isometric extension of $K$.

Let $\bar x$ be a geometric point of $X$. Then one has a functor \[F^{\mbb L}_{\bar
x}:\CovtempL(X)\to\Set\] which maps a $\mbb L$-tempered cover $S\to X$ to the set
$S_{\bar x}$.

The $\mbb L$-tempered fundamental group of $X$ pointed at $\bar x$ is
\[\gtempL(X,\bar x)=\Aut F^{\mbb L}_{\bar x}.\]
When $X$ is a smooth algebraic $K$-variety, $\CovtempL(X^{\an})$
and $\gtempL(X^{\an},\bar x)$ will also be denoted simply by
$\CovtempL(X)$ and $\gtempL(X,\bar x)$.

By considering the stabilizers $(\Stab_{F^{\mbb L}_{\bar x}(S)}(s))_{S\in \CovtempL(X),s\in F^{\mbb L}_{\bar x}(S)}$ as a basis of open subgroups of $\gtempL(X,\bar
x)$, $\gtempL(X,\bar x)$ becomes a topological group. It is a prodiscrete topological group.

When $X$ is algebraic, $K$ of characteristic zero and has only countably
many finite extensions in a fixed algebraic closure $\overline K$,
$\gtempL(X,\bar x)$ has a countable fundamental system of
neighborhood of $1$ and all its discrete quotient groups are finitely
generated~(\cite[prop. III.2.1.7]{andre1}).
When $\mbb L$ is the set of all primes, we often forget it in the
notations.

It should be remarked that in general, for a given $\mbb L$, one
cannot recover $\gtempL(X,\bar x)$ from $\gtemp(X,\bar x)$.  For example, let us consider an Enriques surface $X$ over a nonarchimedean field of residual characteristic zero. One has $\ga(X)=\mbf Z/2\mbf Z$ and $X$ has a unique nontrivial connected finite cover; it is given by a $K3$ surface $Y$. The surfaces $X$ and $Y$ have a semistable reduction, and according to \cite{berk2}, $X^{\an}$ and $Y^{\an}$ are homotopy equivalent to the dual simplicial sets of their semistable reduction. The possible simplicial sets are given by \cite{enriquessurfaces}. For the $K3$ surface $Y$, this dual simplicial set is always simply connected and therefore $\gtemp(X)=\ga(X)=\mbf Z/2\mbf Z$. If $X$ is an Enriques surface with good reduction, $\pi_1^{\emptyset\text{-temp}}(X,\bar x)=\gtop(X,x)=\{1\}$. If the reduction of $X$ is totally degenerate, \ie all the irreducible components of a semistable reduction are projective planes, the dual simplicial set is homotopy equivalent to a real projective plane and $\pi_1^{\emptyset\text{-temp}}(X,\bar x)=\gtop(X,x)=\mbf Z/2\mbf Z$. Therefore, two Enriques surface have isomorphic tempered fundamental groups but can have different $\emptyset$-tempered fundamental groups.

\smallskip

If $\bar x$ and $\bar x'$ are two geometric points, then $F^{\mbb L}_{\bar
x}$ and $F^{\mbb L}_{\bar x'}$ are (non canonically) isomorphic
(\cite[th. 2.9]{dJ1}). Thus, as usual, the tempered fundamental group
depends on the basepoint only up to inner automorphism (this topological
group, considered up to conjugation, will sometimes be denoted simply
$\gtempL(X)$).

The full subcategory of tempered covers $S$ for which $F^{\mbb L}_{\bar
x}(S)$ is $\mbb L$-finite is equivalent to $\Covalg(S)^{\mbb L}$, hence
\[\gtempL(X,\bar x)^{\mbb L}=\ga(X,\bar x)^{\mbb L}\] (where
 $(\ )^{\mbb L}$ denotes the pro-$\mbb L$ completion).

For any morphism $X\to Y$, the pullback defines a functor
$\CovtempL(Y)\to\CovtempL(X)$. If $\bar x$ is a geometric
point of $X$ with image $\bar y$ in $Y$, this gives rise to a continuous
homomorphism \[\gtempL(X,\bar x)\to\gtempL(Y,\bar y)\]
(hence an outer morphism $\gtempL(X)\to\gtempL(Y)$).

One has the analog of the usual Galois correspondence:
\begin{thm}[{\cite[th. III.1.4.5]{andre1}}]  \label{galcorr} $F^{\mbb
L}_{\bar x}$ induces an equivalence of categories between the category of
direct sums of $\mbb L$-tempered covers of $X$ and the category
$\gtempL(X,\bar x)\tSet$ of discrete sets endowed with a
continuous left action of $\gtempL(X,\bar x)$.\end{thm}

If $S$ is a $\mbb L$-finite Galois cover of $X$, its universal
topological cover $S^{\infty}$ is still Galois and every connected $\mbb
L$-tempered cover is dominated by such a Galois $\mbb L$-tempered cover.

If $((S_i,\bar s_i))_{i\in \mbf N}$ is a cofinal projective system (with
morphisms $f_{ij}:S_i\to S_j$ which maps $s_i$ to $s_j$ for $i\geq j$) of
geometrically pointed Galois $\mbb L$-finite étale covers of $(X,\bar
x)$, let $((S^{\infty}_i,\bar s^{\infty}_i))_{i\in\mbf N}$ be the
projective system of its
pointed universal topological covers (the transition maps will be
denoted by $f^{\infty}_{ij}$). It induces a projective system
$(\Gal(S_i^\infty/X))_{i\in\mbf N}$ of discrete groups. For every $i$,
$\Gal(S_i^\infty/X)$ can be identified with $F^{\mbb L}_{\bar
  x}(S_i^{\infty})$: this gives us compatible morphisms $\gtempL(X,\bar
x)\to\Gal(S_i^\infty/X)$.
Then, thanks to~\cite[lem. III.2.1.5]{andre1},
\begin{prop}\label{limproj} \[\gtempL(X,\bar x)\to\varprojlim \Gal(S^{\infty}_i/X)\] is an isomorphism.\end{prop}

In a more categorical way, we have a fibered category
$\Dtop(X)\to\Covalg(X)$, where the fiber $\Dtop(X)_S$ in an algebraic
cover $S$ of $X$ is $\Covtop(X)$.

Since algebraic covers are of effective descent for tempered covers,
the full subcategory of tempered covers $T$ of $X$ such that $T_S\to S$
is a topological cover is naturally equivalent to the category
$\DDtemp_S$ of
descent data in the fibered category $\Dtop(X)$ with respect to $S\to
X$.

If $``\varprojlim"\ S_i$ is a universal procover of $(X,x)$, one gets a
natural equivalence \[\Covtemp(X)=\projLim_i \DDtemp_{S_i}\]
In particular one can recover the tempered fundamental group from the
fibered category $\Dtop(X)\to\Covalg(X)$.

If $S\to S$ is an isomorphism, the induced functor
$\DDtemp_S\to\DDtemp_S$ is naturally isomorphic to the identity. Thus if
$\alpha: ``\varprojlim_i"\ S_i\to ``\varprojlim_i{}"\ S_i$ is an automorphism of
the universal pro-cover, the induced functor $\projLim_i
\DDtemp_{S_i}\to\projLim_i\DDtemp_{S_i}$ is naturally isomorphic to the
identity. Thus the construction does not depend of the choice of the
universal procover.

\smallskip

To give a more stacky and functorial description, let us consider $\Covalg(X)$ with its
canonical topology.

Let $\Dtemp(X)\to\Covalg(X)$ be the fibered category whose fiber over $U$ is
the category $\Covtemp(U)$ of tempered covers of $U$. Then $\Dtemp(X)$
is a stack. The fully faithful cartesian functor of fibered categories $\Dtop(X)\to\Dtemp(X)$ induces a
fully faithful cartesian functor of stacks $\Dtop(X)^a\to\Dtemp(X)$ where $\Dtop(X)^a$ is the stack
associated to $\Dtop(X)$. Since a tempered cover is a topological cover locally on
$\Covalg(X)$, this functor is in fact an equivalence~(\cite[th. II.2.1.3]{giraud}).

In a similar way:
\begin{prop}\label{stackydfn} The stack $(\Dtop(X)|_{\Covalg(X)^{\mbb
    L}})^a$ is the stack $\DtempL(X)$ of $\mbb L$-tempered
covers on $\Covalg(X)^{\mbb L}$.\end{prop}
\subsection{Homotopy types of analytic curves}
Let
 $K$ be a complete nonarchimedean field with separably closed residue field $k$ and let $O_K$ be its ring of integers.
Let $X\to
O_K$ be a proper semistable curve with smooth generic fiber. There is a canonical embedding $|\mbb G(X_k)|\to
X^{\an}_\eta$ which is a homotopy equivalence (\cite[th. 4.3.2]{berk}). If $K'$ is a complete isometric extenstion of $K$ with separably closed field $k'$, then the following diagram is commutative:
\[\xymatrix{|\mbb G(X_{k'})|\ar[d]\ar[r] & X^{\an}_{K'}\ar[d]\\
 |\mbb G(X_k)|\ar[r] & X^{\an}_K}
\]
Moreover, if $U$ is any dense Zariski open subset of $X_{\eta}$, $|\mbb
G(X_k)|$ is mapped into $U^{\an}$ and $|\mbb G(X_k)|\to U^{\an}$ is still a
homotopy equivalence.

If $X\to O_K$ is a semistable log curve and $X'\to X$ is a két morphism
such that $X'$ is still a semistable curve, the following diagram is commutative:
\[\begin{array}{ccc} |\mbb G(X'_k)| & \to & {X'}_{\eta}^{\an}\\ \dar & & \dar\\
|\mbb G(X_k)| & \to & X_{\eta}^{\an}\end{array}\]
\subsection{Cospecialization of tempered fundamental groups}
Let $K$ be a complete discretely valued field. Let $O_K$ be the ring of integers of $K$.
Let $S\to O_K$ be a morphism of fs log schemes of finite type.
Let $S_{\tr}$ be the open locus of $S$  where the log structure is trivial
($S_{\tr}\subset S_\eta$). Let $\fk S$ be the completion of $S$ along its
closed fiber. Then $\fk S_\eta$ is an analytic domain of $S^{\an}$.
Let $S_0=\fk S_\eta\cap S_{\tr}^{\an}\subset S^{\an}$.

Let $\tilde \eta$  be a $K'$-point of $S_0$ where $K'$ is a complete extension of $K$. One has a canonical morphism of log schemes $\Spec O_{K'}\to S$ where $\Spec O_{K'}$ is endowed with the log structure given by $O_{K'}\backslash\{0\}\to O_{K'}$. The \emph{log reduction} $\tilde s$ of $\tilde \eta$ is the log point of $S$
corresponding to the special point of $\Spec O_{K'}$ with the inverse image of the log structure of $\Spec O_{K'}$. If $K'$ has discrete valuation, then $\tilde s$ is a fs log point. If $K'$ is algebraically closed, $\tilde s$ is a geometric log point.

\begin{dfn} The category $\widetilde\Pt^{\an}(S)$ is the category whose objects are the geometric
  points $\bar \eta$ of $Y_0$ such that $\mcal H(\eta)$ is discretely valued (where $\eta$
  is the underlying point of $\bar \eta$) and $\Hom_{\Pt^{\an}(S)}(\bar \eta,\bar
  \eta')$ is the set of két specializations in $S_k$ from the log reduction
  $\bar s$ of $\bar \eta$ to the log reduction $\bar s'$ of $\bar \eta'$ such that there  exists some specialization $\bar \eta\to \bar \eta'$ of geometric points in the sense of algebraic étale topology  for which the following diagram commutes:
\[\xymatrix{ \bar \eta\ar[r]\ar[d] & \bar s\ar[d]\\ 
\bar \eta'\ar[r] & \bar s'
}\]
The category $\Pt^{\an}_0(S)$ is the category obtained from $\widetilde\Pt^{\an}(S)$ by inverting the class of morphisms $\bar \eta\to \bar \eta'$ such that $\overline M_{S,\bar s'}\to\overline M_{S,\bar s}$ is an isomorphism.
\end{dfn}
Let $\OutGptop $ be the category of topological groups with outer morphisms.

\begin{thm}\label{cospcourbes}Let $O_K$ be a complete discretely valued ring of residue characteristic $p\geq 0$, let $\mbb L$ be a set of integers such that $p\notin\mbb L$. Let $S\to \Spec O_K$ be a morphism of fs log schemes of finite type and let $X\to S$ be a proper log
semistable curve. Let $U$ be the open locus of $X$ where the log structure is trivial. Then
there is a functor $\gtempL(U_{(\cdot)}):\Pt^{\an}_0(S)^{\op}\to\OutGptop $
sending $\bar \eta$ to $\gtempL(U_{\bar \eta})$.
\end{thm}
\dem
Let $\bar \eta_2\to\bar \eta_1$ be a morphism of $\widetilde\Pt^{\an}(Y)$. Let us construct a cospecialization morphism $\gtempL(U_{\bar \eta_1})\to\gtempL(U_{\bar \eta_2})$, which is an isomorphism if $\overline M_{S,\bar s_1}\to\overline M_{S,\bar s_2}$ is an isomorphism.

One has a cospecialization functor
\[F:\KCovgeom(X_{s_1}/s_1)^{\mbb L}\to\KCovgeom(X_{s_2}/s_2)^{\mbb L}\] which factors through
$\KCovgeom(X_{T}/T)^{\mbb L}$ where $T$ is the strict localization at $s_1$.\\
The cospecialization functor
$\KCovgeom(X_{s_i}/s_i)^{\mbb L}\to\Covalg(U_{\bar \eta_i})$ is an equivalence since $\eta_i\in \eta_{\tr}$ (\ref{speclogdvrgeom}). If one choses a specialization $\bar \eta_2\to\bar \eta_1$ above $\bar s_2\to\bar s_1$, then one can apply \cite[cor. XIII.2.9]{sga} to $U_K\subset X_K\to S_K$: one gets that the functor
$\Covalg(U_{\bar y_i})^{\mbb L}\to\Covalg(U_{\bar y_2})^{\mbb L}$ is
also an equivalence. Thus $F$ is an
equivalence.

Let $Y_{\bar s_1}$ be a log geometric
két cover of $X_{s_{1}/(s_1,\bar s_1)}$ and let $Y_{\bar s_2}$ (resp. $Y_1$, $Y_2$)
be the corresponding log geometric két cover of $X_{(s_{2},\bar s_2)}$ (resp. $U_{\bar
  \eta_1}$, $U_{\bar \eta_2}$).

\smallskip

There are maps (functorially in $Y$):
\[|Y^{\an}_{\bar s_1}|\leftarrow |\mbb G(Y_{\bar s_1})|\to |\mbb G(Y_{\bar s_2})|\to
|Y^{\an}_{\bar s_2}|\]
where the first and third map are  the embedding of the skeleton of an
anlytic curve. The first and third map are therefore homotopy equivalences.

One thus gets a morphism of homotopy types
$|Y^{\an}_{\bar s_1}|\to|Y^{\an}_{\bar s_2}|$ functorially in $Y$. According to proposition~\ref{topcospket}, if $\overline M_{S,\bar s_1}\to\overline M_{S,\bar s_2}$ is an isomorphism, $|Y^{\an}_{\bar s_1}|\to|Y^{\an}_{\bar s_2}|$ is an isomorphism of homotopy types.

With the notations of proposition~\ref{stackydfn}, one thus gets a functor of fibered categories:
\[\begin{array}{ccc} \Dtop(U_{\bar \eta_2}) & \to & \Dtop(U_{\bar \eta_1})\\
\dar & & \dar\\
 \Covalg(U_{\bar \eta_{2}})^{\mbb L} & \simeq & \Covalg(U_{\bar
   \eta_{1}})^{\mbb L}\end{array}\]
Using proposition~\ref{stackydfn}, this induces a functor of associated stacks:
\[\begin{array}{ccc} \DtempL(U_{\bar \eta_2}) & \to & \DtempL(U_{\bar
    \eta_1})\\
\dar & & \dar\\
 \Covalg(U_{\bar \eta_{2}})^{\mbb L} & \simeq & \Covalg(U_{\bar
    \eta_{1}})^{\mbb L}\end{array}\]
By taking the global sections one gets a functor:
\[\CovtempL(U_{\bar \eta_2})\to\CovtempL(U_{\bar \eta_1}),\]
which is an equivalence if $\overline M_{S,\bar s_1}\to\overline M_{S,\bar s_2}$ is an isomorphism.
It induces a cospecialization outer morphism of tempered fundamental groups
\[\gtempL(U_{\bar \eta_1})\to\gtempL(U_{\bar \eta_2}),\]
which is an isomorphism if $\overline M_{S,\bar s_1}\to\overline M_{S,\bar s_2}$ is an isomorphism.

Let $\bar \eta_3\to\bar\eta_2$ be a morphism of $\widetilde\Pt^{\an}(Y)$.  According to corollary \ref{logsp2}, the diagram
\begin{equation}\label{cospfonct}\xymatrix{\KCovgeom(X_{s_1}/(s_1,\bar s_1))^{\mbb L}\ar[r]^{F_{12}}\ar[rd]^{F_{13}} & \KCovgeom(X_{s_2}/(s_2,\bar s_2))^{\mbb L}\ar[d]^{F_{23}} \\
 & \KCovgeom(X_{s_3}/(s_3,\bar s_3))^{\mbb L}
}
\end{equation}
is 2-commutative. Let $Y_{\bar s_3}$ be the log geometric két cover of $X_{s_3}/(s_3/\bar s_3)$ corresponding to $Y_{\bar s_2}$ and let $Y_{\bar \eta_3}$ be the corresponding cover of $X_{\bar \eta_3}$. The diagram
\[\xymatrix{|\mbb G(Y_{\bar s_1})|\ar[r]\ar[rd] & |\mbb G(Y_{\bar s_2})|\ar[d] \\
 & |\mbb G(Y_{\bar s_3})|
}
\]
is commutative according to \ref{topcospket}, and therefore the diagram of homotopy types
\[\xymatrix{|Y^{\an}_{\bar \eta_1}|\ar[r]\ar[rd] & |Y^{\an}_{\bar s_2}|\ar[d] \\
 & |Y^{\an}_{\bar s_3}|
}
\]
is also commutative. One thus gets a 2-commutative diagram 
\[\xymatrix{\Dtop(U_{\bar \eta_3})\ar[r]\ar[rd] & \Dtop(U_{\bar \eta_2})\ar[d] \\
 & \Dtop(U_{\bar \eta_1})
}
\] of fibered categories above the inverse of \eqref{cospfonct}. By taking the global sections of the associated functor, one gets that the diagram
\[\xymatrix{\CovtempL(U_{\bar \eta_3})\ar[r]\ar[rd] & \CovtempL(U_{\bar \eta_2})\ar[d] \\
 & \CovtempL(U_{\bar \eta_1})
}
\] is 2-commutative, which proves the functoriality of $\gtempL(U_{(\cdot)})$.
\findem

\begin{rem}
Such a functor cannot exist if $p\neq 0$ and $\mbb L$ is the set of all
primes. Consider a moduli space $S$ of stable curves over $\Spec \mbf Z$, endowed with its
canonical log structure. By \cite[Th. 5.1.7]{olsson}, $S$ classifies vertical stable log curves. Let $C\to S$ be the universal log curve. If $\bar s$ is a geometric point of $S$, $\overline M_{S,\bar s}$ can be identified with $\mbf N^I$, where $I$ is the set of double points of $C_{\bar s}$, in such a way that for any compatible local chart $U=\Spec A\stackrel{\phi}{\to} \Spec \mbf Z[\mbf N^I]$ of $S$ at $\bar s$ modeled on $\mbf N^I$, locally on the étale topology around the double point $x$ of $C_{\bar s}$, the universal curve $C\to S$ has a chart
\[\xymatrix{\mbf N^I\oplus\mbf Nu\oplus\mbf Nv/(u+v=e_x) \ar[r] &  A[uv]/(uv-\phi^*(e_x))\\
\mbf N^I \ar[r]^{\phi^*}\ar[u] & A \ar[u]
}
\]
 where $e_x\in \mbf N^I$ is defined by $(e_x)_x=1$ and $(e_x)_{x'}=0$ if $x'\neq x\in I$.
Let $\bar s$ be a geometric point in $S_k$ such that the corresponding stable curve is totally degenerate. Let $I$ be the set of double points of $C_{\bar s}$. The set $I$ is non empty and therefore $\overline M_{S,\bar s}$ is nontrivial.
Let us choose a local chart $\phi:U=\Spec A\to\Spec\mbf Z[\mbf N^I]$ of $S$ at $\bar s$ compatible with the identification of $\overline M_{S,\bar s}$ with $\mbf N^I$ mentioned before.
Choose two different morphisms $a_1,a_2:\mbf N^I\to \mbf N$ such that the preimage of $0$ by $a_1$ and $a_2$ is $0$. Put on $\bar s$ the log structure $k(\bar s)^*\oplus \mbf N$: one gets a log point $s_0$. Given the chart $\phi$, the morphisms $a_1,a_2$ define two morphisms of fs log schemes $s_0\to S$: we denote by $s_1$ and $s_2$ the corresponding fs log points of $S$. By choosing a uniformizer $\pi$ of $O_K$, the morphism of schemes $\bar s\to \Spec O_K$ can be enriched in a morphism of fs log schemes $s_0\to\Spec O_K$ by sending $\pi$ to $(0,1)\in k(\bar s)^*\oplus \mbf N$. Thus $s_1$ and $s_2$ are lifted as fs log points of $S\times_{\Spec \mbf Z}\Spec O_K$.
Let $\eta_1$ and $\eta_2$ be
discretely valued points of $S^{\an}_{K}$ whose log reductions are $s_1$ and $s_2$.
Then, \'etale locally at $x$, $C_{O_{\mcal H(\eta_i)}}$ is isomorphic to $\Spec O_K[u,v]/(uv-\pi^{a_i(e_x)})$. Since $a_i(e_x)>0$ for every $x$, $C_{\eta_i}$ is a smooth curve; since $C_s$ is totally degenerate, $C_{\eta_i}$ is a Mumford curve. The length of the corresponding edge of the graph of $C_{\eta_i}$ is $-\log_p |\pi|^{a_i(e_x)}$. Since $a_1\neq a_2$
the
two Mumford curves $C_{\eta_1}$ and $C_{\eta_2}$ have different metric on the
graph of their stable model, and thus have non isomorphic tempered
fundamental groups (\cite{metricmumford}). But the two geometric log points $s_1$ and $s_2$ are
isomorphic with respect to specialization for két topology since they lie above the same Zariski point.
\end{rem}

\providecommand{\bysame}{\leavevmode\hbox to3em{\hrulefill}\thinspace}
\providecommand{\MR}{\relax\ifhmode\unskip\space\fi MR }
\providecommand{\MRhref}[2]{%
  \href{http://www.ams.org/mathscinet-getitem?mr=#1}{#2}
}
\providecommand{\href}[2]{#2}

\end{document}